
  \magnification 1200
  \input pictex


  \newcount\fontset
  \fontset=1
  \def \dualfont#1#2#3{\font#1=\ifnum\fontset=1 #2\else#3\fi}

  \dualfont\bbfive{bbm5}{cmbx5}
  \dualfont\bbseven{bbm7}{cmbx7}
  \dualfont\bbten{bbm10}{cmbx10}

  \font \eightbf = cmbx8
  \font \eighti = cmmi8 \skewchar \eighti = '177
  \font \eightit = cmti8
  \font \eightrm = cmr8
  \font \eightsl = cmsl8
  \font \eightsy = cmsy8 \skewchar \eightsy = '60
  \font \eighttt = cmtt8 \hyphenchar\eighttt = -1
  \font \msbm = msbm10
  
  \font \sixi = cmmi6 \skewchar \sixi = '177
  \font \sixrm = cmr6
  \font \sixsy = cmsy6 \skewchar \sixsy = '60
  \font \tensc = cmcsc10

  \scriptfont \bffam = \bbseven
  \scriptscriptfont \bffam = \bbfive
  \textfont \bffam = \bbten

  \newskip \ttglue

  \def \eightpoint {\def \rm {\fam0 \eightrm }%
  \textfont0 = \eightrm
  \scriptfont0 = \sixrm \scriptscriptfont0 = \fiverm
  \textfont1 = \eighti
  \scriptfont1 = \sixi \scriptscriptfont1 = \fivei
  \textfont2 = \eightsy
  \scriptfont2 = \sixsy \scriptscriptfont2 = \fivesy
  \textfont3 = \tenex
  \scriptfont3 = \tenex \scriptscriptfont3 = \tenex
  \def \it {\fam \itfam \eightit }%
  \textfont \itfam = \eightit
  \def \sl {\fam \slfam \eightsl }%
  \textfont \slfam = \eightsl
  \def \bf {\fam \bffam \eightbf }%
  \textfont \bffam = \bbseven
  \scriptfont \bffam = \bbfive
  \scriptscriptfont \bffam = \bbfive
  \def \tt {\fam \ttfam \eighttt }%
  \textfont \ttfam = \eighttt
  \tt \ttglue = .5em plus.25em minus.15em
  \normalbaselineskip = 9pt
  \def \MF {{\manual opqr}\-{\manual stuq}}%
  \let \sc = \sixrm
  \let \big = \eightbig
  \setbox \strutbox = \hbox {\vrule height7pt depth2pt width0pt}%
  \normalbaselines \rm }



  \newcount \secno \secno = 0
  \newcount \stno \stno =0
  \newcount \eqcntr \eqcntr=0

  \def \ifn #1{\expandafter \ifx \csname #1\endcsname \relax }

  \def \track #1#2#3{\ifn{#1}\else {\tt\ [#2 \string #3] }\fi}

  \def \laberr#1#2{\message{*** RELABEL CHECKED FALSE for #1 ***}
      RELABEL CHECKED FALSE FOR #1, EXITING.
      \end}

  \def \seqnumbering {\global \advance \stno by 1 \global
    \eqcntr=0 \number \secno .\number \stno }

  \def \current {\number \secno
    \ifnum \number \stno = 0\else .\number \stno \fi }

  \def \eqmark #1 {\global \advance\eqcntr by 1
    \edef\a{\current.\number\eqcntr}
    \eqno {(\a)}
    \syslabel{#1}{\a}
    \track{showlabel}{*}{#1}}

  \def \syslabel#1#2{\global \expandafter \edef \csname
    #1\endcsname {#2}}

  \def \fcite#1#2{\syslabel{#1}{#2}\lcite{#2}}

  \def \label #1 {%
    \ifn {#1}%
      \syslabel{#1}{\current}%
    \else
      \edef\a{\expandafter\csname #1\endcsname}%
      \edef\b{\current}%
      \ifx \a \b \else \laberr{#1=(\a)=(\b)} \fi
      \fi
    \track{showlabel}{*}{#1}}

  \def \lcite #1{(#1\track{showcit}{$\bullet$}{#1})}

  \def \cite #1{[{\bf #1}\track{showref}{\#}{#1}]}

  \def \scite #1#2{{\rm [\bf #1\track{showref}{\#}{#1}{\rm \hskip 0.7pt:\hskip 2pt #2}\rm]}}


 \def \Headlines #1#2{\nopagenumbers
    \advance \voffset by 2\baselineskip
    \advance \vsize by -\voffset
    \headline {\ifnum \pageno = 1 \hfil
    \else \ifodd \pageno \tensc \hfil \lcase {#1} \hfil \folio
    \else \tensc \folio \hfil \lcase {#2} \hfil
    \fi \fi }}

  \def \Date #1 {\footnote {}{\eightit Date: #1.}}


  \def \lcase #1{\edef \auxvar {\lowercase {#1}}\auxvar }

  \def \goodbreak {\vskip0pt plus.1\vsize \penalty -250 \vskip0pt
plus-.1\vsize }

  \def \section #1{\global\def \SectionName{#1}\stno = 0 \global
\advance \secno by 1 \bigskip \bigskip \goodbreak \noindent {\bf
\number \secno .\enspace #1.}\medskip \noindent \ignorespaces}

  \long \def \sysstate #1#2#3{\medbreak \noindent {\bf \seqnumbering
.\enspace #1.\enspace }{#2#3\vskip 0pt}\medbreak }
  \def \state #1 #2\par {\sysstate {#1}{\sl }{#2}}
  \def \definition #1\par {\sysstate {Definition}{\rm }{#1}}
  \def \remark #1\par {\sysstate {Remark}{\rm }{#1}}


  \def \proof {\medbreak \noindent {\it Proof.\enspace }}
  \def \proofend {\ifmmode \eqno \square \else \hfill \square
\looseness = -1 \medbreak \fi }

  \def \$#1{#1 $$$$ #1}
  \def \=#1{\buildrel #1 \over =}

  \def \Item #1{\smallskip \item {#1}}
  \newcount \zitemno \zitemno = 0

  \def \izitem {\zitemno = 0}
  \def \zitemplus {\global \advance \zitemno by 1\relax}
  \def \rzitem{\romannumeral \zitemno}
  \def \rzitemplus {\zitemplus \rzitem}
  \def \zitem {\Item {{\rm(\rzitemplus)}}}

  \newcount \nitemno \nitemno = 0
  
  \def \nitem {\global \advance \nitemno by 1 \Item {{\rm(\number\nitemno)}}}

  \newcount \aitemno \aitemno = 0
  \def\boxlet#1{\hbox to 6.5pt{\hfill #1\hfill}}
  
  \def \aitem {\Item {(\ifcase \aitemno \boxlet a\or \boxlet b\or
\boxlet c\or \boxlet d\or \boxlet e\or \boxlet f\or \boxlet g\or
\boxlet h\or \boxlet i\or \boxlet j\or \boxlet k\or \boxlet l\or
\boxlet m\or \boxlet n\or \boxlet o\or \boxlet p\or \boxlet q\or
\boxlet r\or \boxlet s\or \boxlet t\or \boxlet u\or \boxlet v\or
\boxlet w\or \boxlet x\or \boxlet y\or \boxlet z\else zzz\fi)} \global
\advance \aitemno by 1}

  \newcount \footno \footno = 1
  \newcount \halffootno \footno = 1
  \def \footcntr {\global \advance \footno by 1
  \halffootno =\footno
  \divide \halffootno by 2
  $^{\number\halffootno}$}


  \def \N {{\bf N}}
  
  \def \<{\left \langle \vrule width 0pt depth 0pt height 8pt }
  \def \>{\right \rangle }  
  
  \def \ds{\displaystyle}
  \def \and {\hbox {,\quad and \quad }}
  \def \calcat #1{\,{\vrule height8pt depth4pt}_{\,#1}}

  \def \for #1{,\quad \forall\,#1}
  \def \square {\hbox {$\sqcap \!\!\!\!\sqcup $}}
  \def \rtimes {{\hbox {\msbm o}}}
  \def \stress #1{{\it #1}\/}
  \def \inv {^{-1}}
  \def \*{\otimes}

  \newcount \bibno \bibno =0
  \def \newbib #1{\global\advance\bibno by 1 \edef #1{\number\bibno}}
  \def \bibitem #1#2#3#4{\smallskip \item {[#1]} #2, ``#3'', #4.}
  \def \references {
    \begingroup
    \bigskip \bigskip \goodbreak
    \eightpoint
    \centerline {\tensc References}
    \nobreak \medskip \frenchspacing }


  \def\t{\theta}
  \def\a{\alpha}
  \def\b{\beta}
  \def\End{{\rm End}}
  \def\w{\omega}
  \def\O{\Omega}
  \def\Q{{\bf Q}}
  \def\R{{\bf R}}
  \def\Z{{\bf Z}}
  \def\L{{\cal L}}   \def\L{L}
  \def\bool#1{\big[#1\big]}
  \def\half{^{^{1\kern-0.7pt/2}}}  
  \def\mhalf{^{^{-1\kern-0.7pt/2}}}  
  \def\inf{\wedge} 
  \def\supp{{\rm supp}}
  \def\sup{\vee}   
  \def\lcm{{\rm lcm}}
  \def\Gpd{{\cal G}}
  \def\lprep{v}   \def\lprep{\sigma}
  \def\prep#1{\lprep_{#1}}
  \def\arep{\pi}

  \def\implies{\ \Longrightarrow\ }
  \def\cl#1#2{C^{#1}_{#2}}
  \def\Cl#1#2#3#4{C\,^{#1,#2}_{#3,#4}}
  \def\W#1#2{W_{#1}\!\left(#2\right)}
  \def\admissible{admissible}
  \def\Admissible{Admissible}
  
  \def\anadmissible{an \admissible}


  \def\P{P}
  \def\G{G}
  \def\g{g}
  \def\cy{\omega}
  \def\S{S}


  \newbib\ArzRena
  \newbib\Deaconu
  \newbib\DeaconuTwo
  \newbib\amena
  \newbib\endo
  \newbib\newsgrp
  \newbib\vershik
  \newbib\Hedlund
  \newbib\Ledrappier
  \newbib\Renault
  \newbib\Yeend


  \Headlines
  {Semigroups of local homeomorphisms}
  {R.~Exel and J.~Renault}

  \null\vskip -1cm
  \centerline{\bf SEMIGROUPS OF LOCAL HOMEOMORPHISMS}
  \smallskip
  \centerline{\bf AND INTERACTION GROUPS}
  \footnote{\null}
  {\eightrm 2000 \eightsl Mathematics Subject Classification:
  \eightrm 
  46L55.
  }

  \bigskip
  \centerline{\tensc 
    R.~Exel\footnote{*}{\eightpoint Partially supported by
CNPq and  CAPES/COFECUB.}
  and 
  J.~Renault\footnote{**}{\eightrm Partially supported by
CAPES/COFECUB.}}

  \bigskip
  \Date{23 Aug 2006}

  \midinsert 
  \narrower \narrower
  \eightpoint \noindent
  Given a semigroup of surjective local homeomorphisms on a compact
space $X$ we consider the corresponding semigroup of *-endomorphisms
on $C(X)$ and discuss the possibility of extending it to an
\stress{interaction group}, a concept recently introduced by the first
named author.  We may also define a \stress{transformation groupoid}
whose C*-algebra turns out to be isomorphic to the crossed product
algebra for the interaction group.  Several examples are considered,
including one which gives rise to a slightly different construction
and should be interpreted as being the C*-algebra of a certain
\stress{polymorphism}.
 
  \endinsert

  \section{Introduction}
  Recall from \cite{\newsgrp} that an \stress{interaction group}
consists of a triple $(A,G,V)$, where $A$ is a unital C*-algebra, $G$
is a discrete group, and $V=\{V_g\}_{g\in G}$ is a collection of
positive, unit preserving linear operators on $A$ satisfying
  \izitem
  \zitem $V_1=id_A$,
  \zitem $V_gV_hV_{h\inv} =  V_{gh}V_{h\inv}$,
  \zitem $V_{g\inv}V_gV_h =  V_{g\inv}V_{gh}$,
  \zitem $V_g(ab)=V_g(a)V_g(b)$, for every $a\in A$, and every $b$ in
the range of $V_{g\inv}$.

\medskip\noindent This concept  is a generalization of automorphism groups
designed to deal with dynamical systems in which the transformations
involved are not invertible and perhaps not even single valued.

We give examples of interaction groups arising from actions of
semigroups and identify, in these examples, the crossed product
algebra  \scite{\newsgrp}{6.2} as a groupoid C*-algebra.

Most of our examples are constructed as follows: let $X$ be a compact
topological space, let $G$ be a discrete group, and let $\t$ be a
right action of a given subsemigroup $P\subseteq G$ on $X$ by means of
surjective local homeomorphisms.
  The usual process of dualization, namely:
  $$
  \a_n(f) = f\circ \t_n
  \for n\in P \for f\in C(X),
  $$
  defines an action $\a$ of $P$ on $C(X)$ by endomorphisms.  The first
hurdle that we face is to find an interaction group $V$ of $G$ on $A$
such that $V_n=\a_n$, for all $n\in P$.  

  If such an interaction group is to be found then $V_{n\inv}$ must
necessarily be a transfer operator \scite{\endo}{2.1} for $\a_n$, and
hence we begin our search by looking for  operators $\L_n$ to fulfil
this role.  They must necessarily have the form
  $$
  \L_n(f)\calcat y = \sum_{\theta_n(x)=y} \cy(n,x) f(x)
  \for f\in C(X)
  \for y\in X,
  $$
  where 
  $
  \cy:  P\times X \to \R_+
  $
  is continuous in the second variable, and normalized in the sense
that
  $
  \sum_{\t_n(x)=y} \cy(n,x) =1,
  $
  for every $n\in P$ and  $y\in X$.

Another necessary condition for the existence of the interaction group
$V$ is that $\L_n$ be anti-multiplicative in the semigroup variable
``$n$",
which amounts to $\cy$ satisfying
  $$
  \cy(nm,x) =  \cy(n,x)  \, \cy(m,\t_n(x))
  \for m,n\in P
  \for x\in X.
  $$
  In other words $\cy$ must be a \stress{cocycle} for the semigroup action.
Assuming that $G=P\inv P$, and that $\cy$ satisfies still another
\stress{coherence} condition, we eventually prove that the formula
  $$
  V_g(f)\calcat y = \sum_{\theta_n(x)=y} \cy(n,x) f(\t_m(x))
  ,\quad g=n\inv m
  \for  f\in C(X)
  \for y\in X,
  $$
  does indeed give an interaction group extending $\a$.

  Notice that the above formula amounts to a weighted average of $f$
on $\t_m(\t_n\inv(\{y\}))$, and hence it is a probabilistic version of
the multi-valued map $f\circ \t_m \circ \t_n\inv$, which could also be
written as $f\circ \t_{n\inv m}$, since $\t$ is a right action, except
that $\t_{n\inv m}$ has no meaning under the present hypothesis.

  Motivated by \cite{\Renault}, \cite{\Deaconu} and \cite{\DeaconuTwo}
we define the \stress{transformation groupoid} relative to the
dynamical system $(X,P,\t)$ to be
  $$
  \Gpd = 
  \big\{(x,\g,y)\in X\times \G\times X : \exists\, n,m\in\P,\ 
    \g = n m\inv,\ \t_n(x)=\t_m(y)\big\},
  $$
  under the operations
  $$
  (x,g,y)(y,h,z) = (x,gh,z)
  \and
  (x,g,y)\inv = (y,g\inv,x).
  $$

When all of the favorable conditions are present, in which case 
the interaction group is available, we show
that the crossed product $C(X)\rtimes_V G$ and the groupoid C*-algebra
$C^*(\Gpd)$ are naturally isomorphic.

We then present a series of examples to which our results may be
applied.  Closing our list of examples we look at the case of certain
polymorphisms, a situation that does not precisely fit within our
general framework, but for which we may also consider both the
interaction group and the transformation groupoid, and prove the
isomorphism between the crossed product and the groupoid C*-algebra.

It is interesting to notice that the existence of the interaction
group extending a given semigroup action depends on the existence of
cocycles satisfying suitable \stress{coherence} conditions (see
\fcite{DefineCoherent}{2.7} below).  On the other hand the transformation
groupoid, and therefore also the associated groupoid C*-algebra may be
defined irrespectively of the existence of cocycles.
  Since the coherence condition mentioned above does not seem to be
studied in a systematic way, it might be an interesting project to
determine, from the point of view of Dynamical Systems, whether or not
a semigroup action admitting a coherent cocycle possesses special
properties at all.

Based on the the concept of \stress{cellular automata} we present, in our last
section, an example of a semigroup action of $\N\times\N$ on
Bernoulli's space which fails to  admit a (non-zero) coherent
cocycle.   We believe that a further study of this example might give
some insight on this so far fuzzy 
situation.
  Moreover it seems that this example is at odds with the last
Proposition in \cite{\Deaconu}.

\section{Semigroup actions}
  Given a compact space $X$ we will let $\End(X)$ denote the semigroup
of all surjective local homeomorphisms
  $$
  T: X\to X
  $$
  under the composition law.

  Let $G$ be a group and let $P$ be a subsemigroup (always assumed to
contain $1$) of $G$.  By a \stress{right action} of $P$ on $X$ we shall
mean a map
  $$
  \t : P \to \End(X),
  $$
  such that $\t(1)=id_X$, and
$\t_n\t_m=\t_{mn}$, for all $n$ and $m$ in $P$ (please
notice the order reversal).

Given such a right action $\t$, and 
denoting by $\End(C(X))$ the semigroup of injective unital
*-endomorphisms of the algebra $C(X)$, we will let
  $$
  \a  : P \to \End(C(X))
  $$
  be given by
  $$
  \a_n(f) = f\circ \t_n
  \for n\in P \for f\in C(X).
  $$
  It is immediate to verify that 
$\a_n\a_m=\a_{nm}$, for all $n$ and $m$ in $P$, so that $\a$ becomes a
semigroup action of $P$ on $C(X)$.

We wish to introduce a \stress{transfer operator} (see
\scite{\endo}{2.1}) $\L_n$ for each $\a_n$ so we fix a map
  $$
  \cy:  P\times X \to \R_+,
  $$
  which is continuous in the second variable, and such that
  $$
  \sum_{\t_n(x)=y} \cy(n,x) =1
  \for n\in P \for y\in X.
  \eqno{(\seqnumbering)}   
  \label SumOne
  $$
  Given  $n\in P$  and $f\in C(X)$, let
  $$
  \L_n(f)\calcat y = \sum_{\theta_n(x)=y} \cy(n,x) f(x)
  \for y\in X.
  $$
   It is easy to show that $\L_n(f)\in C(X)$ and that $\L_n$ defines
a transfer operator for $\a_n$, in the sense that $\L_n$
is a positive linear operator on $C(X)$ satisfying $\L_n(1)=1$, and 
  $$
  \L_n\big(f\a_n(g)\big) =   \L_n(f) g
  \for f,g\in C(X).
  $$

  \state Proposition
  \label TransferLaw
  Suppose that 
  $$
  \cy(nm,x) =  \cy(n,x)  \, \cy(m,\t_n(x))
  \for m,n\in P
  \for x\in X.
  $$
  Then $\L_{nm} = \L_m\L_n$.

  \proof
  Given $f\in C(X)$ we have
  $$
  \L_m\L_n(f)\calcat z =
  \sum_{\theta_m(y)=z}  \cy(m,y)
  \sum_{\theta_n(x)=y}  \cy(n,x) f(x) \$=
  \sum_{\theta_m(y)=z} \ 
  \sum_{\theta_n(x)=y}  \cy(m,\t_n(x))\,\cy(n,x) f(x) =
  \sum_{\theta_{nm}(x)=z}  \cy(nm,x) f(x) =
  \L_{nm}(f)\calcat z.
  \proofend
  $$

\definition
  A map
  $
  \cy:  P\times X \to \R_+
  $
  which is continuous in the second variable and which satisfies both
\lcite{\SumOne} and  the condition in \lcite{\TransferLaw} 
will be called a \stress{normalized cocycle}.

From now on we will assume that $\cy$ is a normalized cocycle.

\bigskip 

Given $n\in P$ and $y\in X$ let us denote by $\cl ny$ the set
  $$
  \cl ny = \big\{x\in X: \t_n(x) = \t_n(y) \big\}.
  $$
  One may then show that the range of $\a_n$ consists of the functions
$f\in C(X)$ which are constant on $\cl ny$ for every $y$.  Let
  $$
  E_n = \a_n \circ \L_n,
  $$
  or, in more explicit terms,
  $$
  E_n(f)\calcat y =
  \sum_{x\in \cl ny}
  \cy(n,x)f(x)
  \for f\in C(X)
  \for y\in X.
  $$
  It is an easy exercise to show that $E_n$ is a conditional
expectation from $C(X)$ onto the range of $\a_n$.

\state Lemma
  \label CommutingE
  Let $m,n\in P$.  Then the following are equivalent
  \izitem
  \zitem  $E_n$ and $E_m$ commute.
  \zitem  For every $x,z\in X$ we have that
  $$
  \sum_{y\in \cl mx\cap \cl nz} \cy(n,y)\, \cy(m,x) =
  \sum_{y\in \cl nx\cap \cl mz} \cy(m,y)\, \cy(n,x).
  $$
  
  \proof
  In this proof we will write
  $$
  E_n(f)\calcat y =
  \sum_{x\in X} \bool{\t_n(x)=\t_n(y)}\cy(n,x) f(x),
  $$
  using brackets to mean boolean value and 
  observing that, although the sum is indexed on the infinite set
$X$, only finitely many summands are non-zero.
  We then have for every $z$ in $X$ that
  $$
  E_n E_m (f) \calcat z =
  \sum_{y\in X}\sum_{x\in X} \bool{\t_n(y)=\t_n(z)}
    \bool{\t_m(x)=\t_m(y)}\cy(n,y)\, \cy(m,x) f(x) \$=
  \sum_{x\in X}\left(\sum_{y\in X}\bool{\t_n(y)=\t_n(z)}
    \bool{\t_m(x)=\t_m(y)}\cy(n,y)\, \cy(m,x)\right) f(x) \$=
  \sum_{x\in X}\left(\sum_{y\in \cl mx\cap \cl nz}\cy(n,y)\,
\cy(m,x)\right) f(x).
  $$
  Interchanging $n$ and $m$ we see that 
  $$
  E_m E_n (f) \calcat z =
  \sum_{x\in X}\left(\sum_{y\in \cl nx\cap \cl mz}\cy(m,y)\,
\cy(n,x)\right) f(x),
  $$
  from where the conclusion follows easily.
  \proofend

We would like to further elaborate on condition
\lcite{\CommutingE.ii}.  Given a finite subset $S\subseteq X$, and
$n\in P$, let us denote by
  $$
  \W n{S} = \sum_{y\in S}\cy(n,y).
  \eqno{(\seqnumbering)}
  \label DefineWn
  $$
  Also, given $x,z\in X$, and $n,m\in P$, let
  $$
  \Cl nmxz = \cl nx\cap \cl mz.
  \eqno{(\seqnumbering)}
  \label DefineBigCl
  $$
  Condition \lcite{\CommutingE.ii} may then be rephrased as follows:

\definition 
  \label DefineCoherent
We shall say that a cocycle $\cy$ is \stress{coherent} if for all 
$x,z\in X$, and for all  $m,n\in P$, one has that 
  $$
  \cy(m,x)\ \W n{\Cl mnxz} =   \cy(n,x)\ \W m{\Cl nmxz}.
  $$

We thus arrive at the first point of contact with the theory of
interaction groups.

  \state{Theorem}  
  \label EndosGiveInteraction
  Let
  \izitem
  \zitem $G$ be a group and $X$ be a compact space,
  \zitem $P$ be a subsemigroup of $G$ such that $G=P\inv P$,
  \zitem $\t:P\to \End(X)$ be a right action, and
  \zitem $\cy$ be a normalized coherent cocycle.
  \medskip\noindent Then there exists a unique interaction group
$V=\{V_g\}_{g\in G}$ on $C(X)$, such that $V_n =\a_n$, and
$V_{n\inv}=\L_n,$ for all $n$ in $P$.  Moreover, if $g=n\inv m$, with
$n,m\in P$, then $V_g=\L_n\a_m$, or more explicitly
  $$
  V_g(f)\calcat y = \sum_{\theta_n(x)=y} \cy(n,x) f(\t_m(x))
  \for  f\in C(X)
  \for y\in X.
  $$

  \proof Follows immediately from the results above and
\scite{\newsgrp}{13.3}.  The  formula for $V_g$ above follows from the
fact that $\L_n\a_n=id_{C(X)}$, and hence
  $$
  V_g =
  \L_n\a_n V_{n\inv m} =
  V_{n\inv}V_n V_{n\inv m} =
  V_{n\inv}V_m = \L_n\a_m.
  \proofend
  $$  

\section{The transformation groupoid}
  Given an interaction group, such as one might obtain from 
Theorem \lcite{\EndosGiveInteraction}, we may consider  its crossed
product algebra
  $$
  C(X)\rtimes_V G,
  $$ 
  as introduced by \scite{\newsgrp}{6.2}.  We would like to show that
this algebra coincides with the C*-algebra of a transformation
groupoid naturally constructed from the given dynamical system, which
we introduce in this section.  This is essentially the groupoid
studied in \cite{\Renault}, \cite{\Deaconu}, \cite{\ArzRena} but, as
we shall see, the extension to actions of arbitrary semigroups
requires some further work.

For the time being we will assume that $G$ is a group, $P$ is a
subsemigroup of $G$, and
  $$
  \t: P \to \End(X)
  $$
  is a right action of $P$ on the compact space $X$.

\state Proposition
  \label IntroduceGroupoid
  Suppose that $\P\inv\P\subseteq\P\P\inv$.    Then the set
  $$
  \Gpd = 
  \big\{(x,\g,y)\in X\times \G\times X : \exists\, n,m\in\P,\ 
    \g = n m\inv,\ \t_n(x)=\t_m(y)\big\}
  $$
  is a groupoid under the operations
  $$
  (x,g,y)(y,h,z) = (x,gh,z)
  \and
  (x,g,y)\inv = (y,g\inv,x).
  $$

  \proof
  Let us prove that the operations defined above do indeed give
elements of $\Gpd$.  With respect to the inversion suppose that
$(x,g,y)\in\Gpd$ and let  $n,m\in \P$ be such that
  $\g = n m\inv$
  and
  $\t_n(x)=\t_m(y).$
  Then obviously
  $\g\inv = m n\inv$ 
  and 
  $\t_m(y)=\t_n(x),$
  so $(y,\g\inv,x)\in\Gpd$.

  If moreover $(y,h,z)\in\Gpd$,  let $p,q\in \P$ be such that
  $h = pq\inv$ and $\t_p(y)=\t_q(z).$
  Given that $\P\inv\P\subseteq\P\P\inv$, write $m\inv p\ = u v\inv$,
with $u,v\in \P$.  So
  $$
  \t_{nu}(x) =
  \t_u(\t_n(x)) =
  \t_u(\t_m(y)) =
  \t_{mu}(y) =
  \t_{pv}(y) =
  \t_v(\t_p(y))=
  \t_v(\t_q(z))=
  \t_{qv}(z).
  $$
  In addition 
  $$
  nu (qv)\inv =
  n u v\inv q\inv =
  n m\inv p q\inv  =
  \g h,
  $$
  thus showing that 
$(x,\g h,z)\in\Gpd$.
  We leave it for the reader to check the other groupoid axioms.
  \proofend  

The next result deals with the topological aspects of $\Gpd$.  We
thank for this Trent Yeend, who proved a similar result in
\scite{\Yeend}{3.6}.

\state Proposition
  \label Topology
  Suppose that $\P\inv\P\subseteq\P\P\inv$.    
  For every $n,m\in P$, and for every open sets $A,B\subseteq X$, let
  $$
  \Sigma(n,m,A,B) = \big\{(x,\g,y)\in\Gpd:
    \g = n m\inv,\ \t_n(x)=\t_m(y),\
    x\in A,\ y\in B
  \big\}.
  $$
  Then the collection of all such subsets is a basis for a topology on
$\Gpd$, with respect to which it is a locally compact \'etale
groupoid.

  \proof
  It is obvious that the  $\Sigma(n ,m,A,B)$ cover $\Gpd$.  Next we
must verify that if
  $$
  (x,g,y) \in \Sigma(n_1,m_1,A_1,B_1) \ \cap \ \Sigma(n_2,m_2,A_2,B_2), 
  $$
  then there exists some $\Sigma(n,m,A,B)$ such that  
  $$
  (x,g,y) \in \Sigma(n,m,A,B) \subseteq \Sigma(n_1,m_1,A_1,B_1) \ \cap
\ \Sigma(n_2,m_2,A_2,B_2).
  $$
  By hypothesis there exists $p_1,p_2\in P$ such that
  $
  n_2\inv n_1=p_2p_1\inv.
  $

  For $i=1,2$, let $U_i$ be an open set with
  $$
  \t_{n_i}(x)=\t_{m_i}(y) \in U_i,
  $$
  and such that $\t_{p_i}$ is injective on $U_i$.  Defining
  \medskip\item{$\bullet$} 
    $ n=n_1p_1$, \ $m=m_1p_1, $
  \medskip\item{$\bullet$} 
    $A = A_1\cap A_2 \cap \t_{n_1}\inv (U_1) \cap \t_{n_2}\inv (U_2),$
  \medskip\item{$\bullet$}   
   $B= B_1\cap B_2 \cap \t_{m_1}\inv (U_1) \cap \t_{m_2}\inv (U_2),$
  \medskip\noindent we claim that $\Sigma(n,m,A,B)$ has the desired
properties.  To see this first observe that 
$g=n_1m_1\inv=n_2m_2\inv$, whence
  $$
  m=m_1p_1=
  m_2n_2\inv n_1p_1 =   m_2p_2p_1\inv p_1 = m_2p_2,
  $$
  and likewise
  $$
  n=n_1p_1=n_2p_2.
  $$
  Thus if $(x',g,y')\in \Sigma(n,m,A,B)$ we have for $i=1,2$
that
  $$
  \t_{p_i}(\t_{n_i}(x')) =
  \t_{n_ip_i}(x') =
  \t_n(x') =  \t_m(y') =
  \t_{m_ip_i}(y') =
  \t_{p_i}(\t_{m_i}(y')).
  $$
  Moreover notice that
  $$
  x'\in A \subseteq  \t_{n_i}\inv (U_i) \implies \t_{n_i}(x')\in U_i,
  $$$$
  y'\in B \subseteq  \t_{m_i}\inv (U_i) \implies \t_{m_i}(y')\in U_i.
  $$
  Since   $\t_{p_i}$ is injective on $U_i$ we have that 
  $\t_{n_i}(x') =\t_{m_i}(y')$, which implies that
  $
  (x',g,y')\in\Sigma(n_i,m_i,A_i,B_i),
  $
  thus showing  that 
  $$
  \Sigma(n,m,A,B) \subseteq 
  \Sigma(n_1,m_1,A_1,B_1) \ \cap \ \Sigma(n_2,m_2,A_2,B_2).
  $$
  In order to see that $(x,g,y)\in\Sigma(n,m,A,B)$, notice that
  $$
  \t_n(x) =
  \t_{n_1p_1}(x) =
  \t_{p_1}(\t_{n_1}(x)) =
  \t_{p_1}(\t_{m_1}(y)) =
  \t_{m_1p_1}(y) =
  \t_m(y).
  $$
  It is also easy to see that $x\in A$ and $y\in B$, so indeed
$(x,g,y)\in\Sigma(n,m,A,B)$.

Before proving that $\Gpd$ is locally compact let for every $n,m\in P$
  $$
  E(n,m) = \big\{(x,y)\in X\times X: \t_n(x) = \t_m(y)\big\},
  $$
  considered as a topological subspace of $X\times X$.  We claim that
the map
  $$
  \iota : (x,y)\in E(n,m)  \mapsto (x,nm\inv,y)\in\Gpd
  $$  
  is continuous.
  To prove that $\iota$ is continuous at a point $(x_0,y_0)\in
E(n,m)$, let $W$ be a neighborhood of $(x_0,nm\inv,y_0)$ in $\Gpd$.
Then there exists a basic open set, say $\Sigma(k,l,A,B)$, such that
  $$
  (x_0,nm\inv,y_0)\in \Sigma(k,l,A,B) \subseteq W.
  $$
  One then must have $nm\inv=kl\inv$.
  Choose $p,q\in P$ such that $k\inv n =
pq\inv$, and let $W$ be an open neighborhood of $\t_k(x_0)=\t_l(y_0)$,
such that $\t_p$ is injective on $W$.    Letting 
  \medskip\item{$\bullet$} $U=A\cap\t_k\inv(W)$,
  \medskip\item{$\bullet$} $V=B\cap \t_l\inv(W)$, and
  \medskip\item{$\bullet$} $Z=(U\times V)\cap E(n,m)$,
  \medskip\noindent
  observe that $Z$ is an open subset of
$E(n,m)$ containing $(x_0,y_0)$.   We claim that 
  $
  \iota(Z)\subseteq W.
  $
  In fact, if $(x,y)\in Z$ then 
  $$
  \t_p(\t_k(x)) =
  \t_{kp}(x) =
  \t_{nq}(x) =
  \t_q(\t_n(x)) =
  \t_q(\t_m(y)) =
  \t_{mq}(y) = \ldots
  $$
  Observe that
  $
  mq = lk\inv n q =  lpq\inv q = lp,
  $
  so the above equals
  $$
  \ldots =
  \t_{lp}(y) =
  \t_p(\t_l(y)).
  $$
  Summarizing we have that $\t_p(\t_k(x)) =  \t_p(\t_l(y))$, but since 
  $$
  x\in U \subseteq  \t_k\inv (W) \implies \t_k(x)\in W,
  $$$$
  y\in V \subseteq  \t_l\inv (W) \implies \t_l(y)\in W,
  $$
  and since $\t_p$ is injective on $W$, we have that
$\t_k(x)=\t_l(y)$, from where one concludes that
  $$
  \iota(x,y) = (x,nm\inv,y) = (x,kl\inv,y) \in \Sigma(k,l,A,B)
\subseteq W,
  $$
  thus proving that $\iota(Z)\subseteq W$, and hence that $\iota$ is
continuous.

Viewing   $\iota$ as a map
  $$
  \iota : E(n,m) \to \Sigma(n,m,X,X) =: \Sigma(n,m),
  $$
  we thus see that $\iota$ is a homeomorphism onto $\Sigma(n,m)$,
because $E(n,m)$ is compact and $\Sigma(n,m)$ is Hausdorff.  Since
$\Sigma(n,m)$ is open by definition we deduce that $\Gpd$ is locally
compact.
  We leave the verification of the remaining properties to the reader.
  \proofend

\section{A semigroup of isometries} 
  Our major goal is to prove the isomorphism between
$C(X)\rtimes_V G$ and $C^*(\Gpd)$.  In order for both of these algebras to
be defined we may invoke \lcite{\EndosGiveInteraction} and \lcite{\Topology}
and hence we must restrict ourselves to a situation in which all of the
relevant hypothesis are satisfied.  We therefore  suppose throughout
that we are under the following:

  \sysstate{Standing Hypotheses}{\rm}{\label StandingHyp
  \izitem
  \zitem $G$ is a discrete group and $P$ is a subsemigroup such that
$1\in P$, and $G=P\inv P = PP\inv$,
  \zitem $\t: P \to \End(X)$ is right action of $P$ on
the compact space $X$, and
  \zitem $\cy$ is a normalized coherent cocycle which never vanishes.}

The proof of the aforementioned isomorphism will be based on the
construction of a suitable covariant representation $(\arep,\lprep)$
\scite{\newsgrp}{4.1} of our interaction group.  In this section we
construct a semigroup of isometries which will later be used to
construct the partial representation $\lprep$.

\state Proposition
  \label SemigroupIsometries
For each $n\in P$, let $\S_n$ be the
element of $C^*(\Gpd)$ given by
  $$
  \S_n(x,\g,y) = 
  \cy(n,x)\half \bool{\g=n}\bool{\t_n(x)=y},
  $$
  where the brackets correspond to boolean value.  Then
$\{\S_n\}_{n\in P}$ is a semigroup of isometries in $C^*(\Gpd)$.

  \proof
  In order to verify that $\S_n$ is an isometry for every $n$ first
observe that
  $$
  \matrix{
  \S_n^*(x,\g,y) & = &  \overline{ \S_n(y,\g\inv,x)}\ = \hfill\cr\cr
  & = &
  \cy(n,y)\half \bool{\g=n\inv}\bool{x=\t_n(y)}.}
  $$
  Therefore 
  $$
  \S_n^*\S_n\calcat{(x,\g,y)} =
  \sum_{(x,h,z)\in\Gpd} \S_n^*(x,h,z) \S_n(z,h\inv\g,y) \$=
  \sum_{(x,h,z)\in\Gpd} 
    \cy(n,z)\half \bool{h=n\inv}\bool{x=\t_n(z)}
    \,\cy(n,z)\half \bool{h\inv\g=n} \bool{\t_n(z)=y} \$=
  \sum_{\t_n(z)=x} 
    \cy(n,z) 
    \bool{\g=1} 
    \bool{\t_n(z)=y} \$=
  \bool{\g=1} \bool{x=y} \sum_{\t_n(z)=x} \cy(n,z) =
  \bool{\g=1} \bool{x=y},
  $$
  thus showing that $\S_n^*\S_n=1$.  We next show that
$\S_n\S_m=\S_{nm}$, for every $n,m\in\P$.
  $$
  \S_n\S_m\calcat{(x,\g,y)} =  
  \sum_{(x,h,z)\in\Gpd} \S_n(x,h,z) \S_m(z,h\inv\g,y) \$=
  \sum_{(x,h,z)\in\Gpd} 
    \cy(n,x)\half\bool{h=n}\bool{\t_n(x)=z}
    \cy(m,z)\half\bool{h\inv\g=m}\bool{\t_m(z)=y} \$=
  \cy(n,x)\half
  \cy(m,\t_n(x))\half\bool{n\inv\g=m}\bool{\t_m(\t_n(x))=y} \$=
  \cy(nm,x)\half
  \bool{\g=nm}\bool{\t_{nm}(x)=y} =
  \S_{nm}(x,\g,y).
  \proofend
  $$

Observe that we have not used the fact that $\cy$ is coherent in the
above proof.

Recall that, whenever $\Gpd$ is an \'etale groupoid, with unit space
$\Gpd^{(0)}$, the algebra $C_0(\Gpd^{(0)})$ sits naturally as a
subalgebra of $C^*(\Gpd)$.  In our case $\Gpd^{(0)}=X$, so we will
henceforth identify $C(X)$ with the corresponding subalgebra of $C^*(\Gpd)$.

We now wish to show that $C^*(\Gpd)$ is generated, as a C*-algebra, by
  $$
  C(X)\cup \{\S_n: n\in P\}.
  $$
  In preparation for this we will
occasionally consider the pointwise (as opposed to the convolution)
product for functions on $\Gpd$.  In order to avoid confusion we will
denote pointwise product by $f\cdot g$, keeping the usual
juxtaposition notation for the convolution product.

In what follows we will say that a subset $E\subseteq\Gpd$ is an
\stress{$r$-section} when $r$ (the groupoid \stress{range} map) is
injective on $E$.

\state Lemma
  \label TwoProducts
  Let $f$ and $g$ be continuous complex valued functions on $\Gpd$,
and suppose that $\supp(f)$ is a compact $r$-section.  Then there
exists $u\in C(X)$ such that
  $$
  g\cdot f = uf.
  $$

  \proof
  Let $K=\supp(f)$, and let $r\inv$ refer to the inverse of the
restriction
  $$
  r: K\to r(K).
  $$
  Since $r$ is continuous and $K$ is compact we have that $r\inv$ is
continuous.
  Define $u_0$ on $r(K)$ by $u_0=g\circ r\inv$ and use Tietze's
Theorem to extend $u_0$ to a continuous function $u$ on $X$.
  Then, for every $\gamma\in\Gpd$ we have
  $$
  (uf)(\gamma) = u(r(\gamma))f(\gamma) = g(\gamma) f(\gamma) = (g\cdot
f)(\gamma).
  \proofend
  $$

The following elementary result will be useful later and is included
for completeness.

\state Lemma
  \label StoneWeierstrass
  Let $\Omega$ be a locally compact space and let ${\cal F}\subseteq
C_c(\Omega)$ be such that for every $x\in\Omega$ there is an $f\in{\cal
F}$ such that $f(x)\neq 0$.  Then
  $$
  C_c(\Omega) = {\sl span}\big\{g f: g\in C_c(\Omega),\ f\in{\cal
F}\big\}.
  $$

  \proof
  Given $h\in C_c(\Omega)$ let $K=\supp(h)$.  For each $x\in K$ choose
$f_x\in{\cal F}$ such that $f_x(x)\neq 0$.  Setting $U_x=\{y\in\Omega:
f_x(y)\neq0\}$ we have that  $\{U_x\}_{x\in K}$ is an open cover of $K$, so
  $$
  K\subseteq U_{x_1}\cup\ldots U_{x_n},
  $$
  for suitable $x_1,\ldots,x_n\in K$.  For simplicity of notation set
$f_i:=f_{x_i}$ and observe that
  $
  \sum_{i=1}^n |f_i|^2 >0
  $
  on $K$.  By Tietze's Theorem choose $g\in C(\Omega)$ such that
  $$
  g(x) = \left(\sum_{i=1}^n |f_i(x)|^2\right)\inv
  \for x\in K.
  $$
  It follows that
  $
  \ds\sum_{i=1}^n (hg\bar f_i) f_i = h.
  $
  \proofend

We thus arrive at a main result:

\state Proposition
  \label SpanGeneratesGpd
  $
  C_c(\Gpd) =
  {\sl span}\big\{u \S_n\S_m^* v: n,m\in \P,\ u,v\in C(X)\big\}.
  $

\proof
  Let ${\cal F}\subseteq C_c(\Gpd)$ be the collection of all functions
of the form $\S_n\S_m^* v$, where $n,m\in\P$, and $v\in C(X)$ is such
that $\t_m$ is injective on $\supp(v)$.

Observe that the support of $\S_n\S_m^* v$, is contained in
$\Sigma(n,m,X,\supp(v))$. In fact, if
  $$
  0 \neq \S_n\S_m^* v\calcat {(x,g,y)} =
  v(y) \cy(n,x)\half\cy(m,y)\half \bool{g=nm\inv}
\bool{\t_n(x)=\t_m(y)},
  $$
  then
  $y\in \supp(v)$, $g=nm\inv$, and $\t_n(x)=\t_m(y)$, which says that
$(x,g,y)$ belongs to the indicated set.

Notice that $\Sigma(n,m,A,B)$ is an $r$-section whenever $\t_m$ is
injective on $B$.  In fact, if
  $$
  (x,nm\inv,y), (x',nm\inv,y') \in\Sigma(n,m,A,B),
  $$
  and $x=x'$, then
  $$
  \t_m(y) = \t_n(x) = \t_n(x') = \t_m(y'),
  $$
  so that $y=y'$ because $\t_m$ is injective on $B$.  It follows that
the support of every   $f\in \cal F$ is an $r$-section.
 
Given that ${\cal F}$ satisfies the hypothesis of
\lcite{\StoneWeierstrass}, as one may easily verify, we conclude that
any $h\in C_c(\Gpd)$ may be written as
  $$
  h= \sum_{i=1}^n g_i\cdot f_i,
  $$
  where $g_i\in C_c(\Gpd)$, and $f_i\in{\cal F}$.  But since $f_i$ is
supported on an $r$-section, as seen above, we have that
  $g_i\cdot f_i = u_if_i$, for some $u_i\in C(X)$, by
\lcite{\TwoProducts}.  Writing $f_i = \S_{n_i}\S_{m_i}^* v_i$, we have
that
  $$
  h=
  \sum_{i=1}^n u_i\S_{n_i}\S_{m_i}^* v_i,
  $$
  as required.
  \proofend

\section{A partial representation}
  We next wish to study the possibility of extending $\S$ to a partial
representation $\lprep$ of $G$ in $C^*(\Gpd)$ which will be part of
the covariant representation we are looking for.  The following
abstract result will suit us well:

\state Proposition
  \label GeneralPiso
  Let $P$ be a subsemigroup of a group $G$ such that $G=\P\inv P$, and
let $\S:P\to B$ be a semigroup of isometries in a C*-algebra $B$.  If
$\S_n\S_n^*$ and $\S_m\S_m^*$ commute for all $n$ and $m$ in $P$ then
there exists a unique *-partial representation
  $$
  \lprep: G \to B
  $$
  such that $\prep{n}=\S_n$, for every $n\in P$.  Moreover, if $g=n\inv
m$, with $n,m\in P$, then   $\prep{g} = \S_n^*\S_m.$

  \proof
  Supposing that $\lprep$ exists let us first address the very last
part of the statement and hence also uniqueness.  Given $g=n\inv m$,
observe that $\prep {m\inv} \prep m = \S_m^* \S_m = 1$, so
  $$
  \prep{g} = 
  \prep{n\inv m} \prep {m\inv} \prep m =
  \prep{n\inv} \prep m =  \S_n^*\S_m.
  $$

  With respect to existence, for every  $g\in G$,  write $g=n\inv m$,
with $n,m\in P$, and set
  $$
  \prep{g} = \S_n^*\S_m.
  $$  
  Let us prove that this is well defined: if $g$ can
also be written as $g=p\inv q$, we claim that there exists $x,y,u,v\in
P$ such that
  $$
  \left\{
  \matrix{
    xp = u, \cr 
    xq = v, \cr
    yn = u, \cr
    ym = v.}\right.
  $$
  In fact, use the hypothesis to write $np\inv=y\inv x$, with $x,y\in
P$, and set $u=xp$, and $v=xq$.  Then the first two equations hold
true and one clearly has that $yn=xp=u$, while
  $$
  ym =
  ynn\inv m =
  xpg=
  xpp\inv q =
  xq = v.
  $$
  One therefore has that
  $$
  \S^*_n\S_m = \S^*_n \S_y^* \S_y \S_m =
  \S_{yn}^*\S_{ym} =
  \S_u^*\S_v =
  \S_{xp}^*\S_{xq} =
  \S_p^*\S_x^*\S_x\S_q =
  \S_p^*\S_q.
  $$
  This proves that $\lprep$ is well defined.
  Let us now show the partial group law, that is,
  $$
  \prep{g}\prep{h}\prep{h\inv} = \prep{gh}\prep{h\inv},
  $$
  for all $g,h\in G$.  Write $g=n\inv m$, and $h=p\inv q$, with
$n,m,p,q\in\P$.  Pick $u,v\in\P$ such that $pm\inv=v\inv u$ and notice
that
  $$
  um=vp.
  $$
  Replacing
  $(n,m)$ by $(un,um)$, and $(p,q)$ by $(vp,vq)$, we may then assume that
$m=p$. We then have
  $$
  \prep{g}\prep{h}\prep{h\inv} = 
  \S^*_n\S_m\S^*_m \S_q\S^*_q\S_m =
  \S^*_n\S_q\S^*_q \S_m\S^*_m\S_m =
  \S^*_n\S_q\S^*_q\S_m =
  \prep{gh}\prep{h\inv}.
  \proofend
  $$

  We now wish to verify that the hypothesis of the above abstract
result does indeed apply in our situation.  Our next Lemma is
designed to isolate the more technical aspects of that verification.
In it we will make extensive use of the notation introduced in
\lcite{\DefineWn} and \lcite{\DefineBigCl}.

\state Lemma
  \label GreatTechnicalComput
  Assuming \lcite{\StandingHyp} let $x,z\in X$.  Then the expressions
  \izitem
  \zitem 
  $
  \W m{\Cl nmxx}\ \W n{\Cl mnxz}, 
  $
  and
  \zitem
  $
  \ds
  \sum_{y\in\Cl mnxz} 
    \cy(m,y)\half \cy(n,y)\half 
    \cy(m,x)\half \cy(n,z)\half
  $
  \medskip\noindent
  are invariant under exchanging the variables $n$ and $m$ in $P$.

\proof With respect to (i) we have to prove that
  $$
 \W m{\Cl nmxx}\ \W n{\Cl mnxz} =  \W n{\Cl nmxx}\ \W m{\Cl nmxz}.
  $$
  For every $y$ we have by \lcite{\DefineCoherent} that
  $$
  \cy(m,y)\ \W n{\Cl mnyz} =   \cy(n,y)\ \W m{\Cl nmyz}.
  $$
  If $y\in\Cl nmxx$ then $\Cl mnyz = \Cl mnxz$, and $\Cl nmyz = \Cl
nmxz$, so the equation displayed above   becomes
  $$
  \cy(m,y)\ \W n{\Cl mnxz} =   \cy(n,y)\ \W m{\Cl nmxz}.
  $$
  Adding up  both sides over all $y\in\Cl nmxx$ we obtain (i).

Let us now deal with (ii).
  If $\Cl mnxz=\emptyset$, then by \lcite{\DefineCoherent} we have that
$\Cl nmxz=\emptyset$ as well (because $\cy$ never vanishes) and hence
the expression in (ii) vanishes regardless of the order in which we
take $n$ and $m$.  On the other hand, given $y\in \Cl mnxz$, we have
that
  $\Cl nmxz\neq\emptyset$, again by
\lcite{\DefineCoherent}, so let $u\in\Cl nmxz$.

Using 
\lcite{\DefineCoherent} we have that
  $$
  \cy(m,y) \W n{\Cl mnyu} = 
  \cy(n,y) \W m{\Cl nmyu}.
  $$
  Notice that
  $
  \Cl mnyu = \Cl mnxx,
  $ and $
  \Cl nmyu = \Cl nmzz,
  $
  so  the above becomes 
  $$
  \cy(m,y) \W n{\Cl mnxx} = 
  \cy(n,y) \W m{\Cl nmzz}.
  \eqmark NewEqtn
  $$
  Working with just the part of the expression in (ii) that
involves the variable $y$, we have
  $$
  \sum_{y\in\Cl mnxz} \cy(m,y)\half \cy(n,y)\half \$=
  \sum_{y\in\Cl mnxz} \cy(n,y)\half \W m{\Cl nmzz}\half 
    \W n{\Cl mnxx}\mhalf \cy(n,y)\half \$=
  \W m{\Cl nmzz}\half \W n{\Cl mnxx}\mhalf 
    \sum_{y\in\Cl mnxz} \cy(n,y) \$=
  \W m{\Cl nmzz}\half \W n{\Cl mnxx}\mhalf 
    \W n{\Cl mnxz}.
  $$ 
  We thus conclude that  (ii) equals
  $$
  \cy(m,x)\half \cy(n,z)\half
  \W m{\Cl nmzz}\half \W n{\Cl mnxx}\mhalf \underline{\W n{\Cl mnxz}}  \$= 
  \cy(m,x)\half \underline{\W n{\Cl mnxz}}\half \cy(n,z)\half
  \W m{\Cl nmzz}\half \W n{\Cl mnxx}\mhalf \underline{\W n{\Cl mnxz}}\half.
  \eqmark AlmostSymmetric
  $$
  Using \lcite{\DefineCoherent}, as it stands, we deduce that
\lcite{\AlmostSymmetric} equals
  $$
  \cy(n,x)\half \W m{\Cl nmxz}\half \cy(n,z)\half \W m{\Cl nmzz}\half 
    \W n{\Cl mnxx}\mhalf \W n{\Cl mnxz}\half = \ldots
  $$  
  Using \lcite{\DefineCoherent} again, this time in the form
  $$
  \cy(n,z)\ \W m{\Cl nmzz} =  \cy(m,z)\ \W n{\Cl nmzz},
  $$
  we conclude that \lcite{\AlmostSymmetric} equals
  $$
  \cy(n,x)\half \W m{\Cl nmxz}\half
  \cy(m,z)\half \W n{\Cl nmzz}\half
  \W n{\Cl mnxx}\mhalf \W n{\Cl mnxz}\half.
  $$
  We want to show this to equal \lcite{\AlmostSymmetric} with
$n$ and $m$ exchanged, namely
  $$
  \cy(n,x)\half \cy(m,z)\half
  \W n{\Cl mnzz}\half \W m{\Cl nmxx}\mhalf \W m{\Cl nmxz},
  $$
  which turns out to be equivalent to verifying that 
  $$
  \W m{\Cl nmxz}\half
   \W n{\Cl mnxx}\mhalf \W n{\Cl mnxz}\half
  =
   \W m{\Cl nmxx}\mhalf \W m{\Cl nmxz},
  $$
  which happens to be  just a rewriting of (i).
  \proofend

  The more involving technical aspects taken care of, we may prove:

  \state Proposition
  \label GotPrep
  Under \lcite{\StandingHyp} there exists a unique partial
representation
  $$
  \lprep : G \to C^*(\Gpd)
  $$
  such that $\prep{n}=\S_n$, for every $n\in P$, where $\{S_n\}_{n\in
P}$ is the semigroup of isometries in $C^*(\Gpd)$ given by
\lcite{\SemigroupIsometries}.
Moreover, if $g=n\inv
m$, with $n,m\in P$, then   $\prep{g} = \S_n^*\S_m.$

\proof
  In view of \lcite{\GeneralPiso} it is enough to verify that
  $\S_n\S_n^*$ and $\S_m\S_m^*$ commute for all $n$ and $m$ in $P$.
Observe that, for all $(x,g,y)\in\Gpd$, we have
  $$
  \S_n\S_n^*(x,g,y) =
  \sum_{(x,h,z)\in\Gpd} \S_n(x,h,z) \S_n^*(z,h\inv g,y) \$=
  \sum_{(x,h,z)\in\Gpd}
    \cy(n,x)\half\bool{h=n}\bool{\t_n(x)=z} 
    \cy(n,y)\half\bool{h\inv g=n\inv}\bool{\t_n(y)=z} \$= 
  \cy(n,x)\half\cy(n,y)\half\bool{g=1}\bool{\t_n(x)=\t_n(y)}.
  \eqmark FormulaSSStar
  $$
  Thus
  $$
  \S_m\S_m^* \S_n\S_n^*(x,g,z) =
  \sum_{(x,h,y)\in\Gpd} \S_m\S_m^*(x,h,y) \S_n\S_n^*(y,h\inv g,z) \$=
  \sum_{(x,h,y)\in\Gpd} 
    \cy(m,x)\half\cy(m,y)\half\bool{h=1}\bool{\t_m(x)=\t_m(y)}
  \kern 4cm $$ \vskip-20pt $$ \kern 5cm 
    \cy(n,y)\half\cy(n,z)\half\bool{h\inv g=1}\bool{\t_n(y)=\t_n(z)}\$=
  \bool{g=1}  \sum_{y\in\Cl mnxz} 
    \cy(m,x)\half\cy(m,y)\half
    \cy(n,y)\half\cy(n,z)\half.
  $$
  By \lcite{\GreatTechnicalComput} we have that this is symmetric in
$m$ and $n$, hence $\S_m\S_m^* \S_n\S_n^*
=\S_n\S_n^*\S_m\S_m^*$.  The result then follows from \lcite{\GeneralPiso}.
\proofend

\section{The isomorphism}
  We shall now employ the conclusions of \lcite{\GotPrep} in order to
construct a covariant representation of the interaction group
$\big(C(X),G,V\big)$ in $C^*(\Gpd)$, which will lead to an isomorphism
between the crossed product $C(X)\rtimes_V G$ and $C^*(\Gpd)$.
As always, we keep \lcite{\StandingHyp} in force.

Consider the canonical representation 
  $$
  \arep: C(X)\to C^*(\Gpd)
  $$
  given
by
  $$
  \arep(f)(x,g,y) = f(x) \bool{g=1}\bool{x=y}
  \for f\in C(X)
  \for (x,g,y)\in\Gpd.
  $$
  
\state Proposition
  \label TheCovarRep
  The pair $(\arep,\lprep)$ is a covariant representation of the
interaction group $\big(C(X),G,V\big)$ in $C^*(\Gpd)$.

\proof
  We are required to show that
  $$
  \prep{g}\arep(f)\prep{g\inv} = \arep(V_g(f)) \prep{g}\prep{g\inv}
  \for f\in C(X)
  \for g\in G.
  \eqmark CovarCond
  $$

\medskip\noindent  {\tensc Case 1:}
  Suppose that $g=n\in P$.  Then, for every $(x,g,y)\in \Gpd$, we have
  $$
  \prep n \arep(f) \prep{n\inv}\calcat{(x,g,y)} =
  \sum_{(x,h,z)\in \Gpd} S_n(x,h,z) f(z) S_n^*(z,h\inv g,y) \$=
  \sum_{(x,h,z)\in \Gpd} \cy(n,x)\half\bool{h=n}\bool{\t_n(x)=z} f(z) 
  \cy(n,y)\half\bool{h\inv g = n\inv}\bool{\t_n(y)=z} \$=
  \cy(n,x)\half
  f(\t_n(x)) 
  \cy(n,y)\half\bool{g = 1}\bool{\t_n(y)=\t_n(x)} \$=
  f(\t_n(x)) \big(S_nS_n^*\big)\calcat{(x,g,y)} =
  \big(\arep(V_n(f))\prep n \prep {n\inv}\big)\calcat{(x,g,y)},
  $$
  thus proving \lcite{\CovarCond}.

\medskip\noindent  {\tensc Case 2:}
  Suppose that $g=n\inv$, where $n\in P$.  Then, for every
$(x,g,y)\in \Gpd$, we have
  $$
  \prep {n\inv} \arep(f) \prep{n}\calcat{(x,g,y)} =
  \sum_{(x,h,z)\in \Gpd} S_n^*(x,h,z) f(z) S_n(z,h\inv g,y) \$=
  \sum_{(x,h,z)\in \Gpd} \cy(n,z)\half\bool{h=n\inv}\bool{\t_n(z)=x}
f(z) 
  \cy(n,z)\half  \bool{h\inv g = n}  \bool{\t_n(z)=y} \$=
  \bool{g=1} \bool{x=y}
  \sum_{z\in\cl nx} \cy(n,z) f(z) =
  \bool{g=1} \bool{x=y}
  \L_n(f)\calcat x,
  $$
  so 
  $$
  \prep {n\inv} \arep(f) \prep{n} =
  \arep(\L_n(f)) = 
  \arep(V_{n\inv}(f)) \prep {n\inv} \prep{n},
  $$
  because $\prep n$ is an isometry.

\medskip\noindent {\tensc Case 3:} For a general $g\in G$, write
$g=n\inv m$, where $n,m\in P$.   By \lcite{\GotPrep} we have 
  $\prep{g}=\prep{n\inv}\prep{m}$, while 
\lcite{\EndosGiveInteraction} gives   $V_g=V_{n\inv}V_m.$
  Therefore
  $$
  \prep{g}\arep(f)\prep{g\inv} = 
  \prep{n\inv}\prep m\arep(f)\prep{m\inv}\prep n =
  \prep{n\inv}\arep(V_m(f))\prep m\prep{m\inv}\prep n = \ldots
  $$
  Observe that
  $$
  \prep m\prep{m\inv}\prep n =
  \prep m\prep{m\inv}\prep n \prep {n\inv}\prep n =
  \prep n \prep {n\inv} \prep m\prep{m\inv} \prep n,
  $$
  so the above equals
  $$
  \ldots = 
  \prep{n\inv}\arep(V_m(f)) \prep n \prep {n\inv} \prep m\prep{m\inv}
\prep n =
  \arep(V_{n\inv}(V_m(f))) \prep{n\inv} \prep n \prep {n\inv} \prep
m\prep{m\inv} \prep n \$=
  \arep(V_g(f)) \prep {n\inv} \prep m\prep{m\inv} \prep n =
  \arep(V_g(f)) \prep g\prep{g\inv}.
  \proofend
  $$

Given the  covariant representation $(\arep,\lprep)$ above we may
use \scite{\newsgrp}{5.3} to 
define a *-homomor\-phism  
  $$
  \widehat{\arep\times\lprep} : {\cal T} \big(C(X),G,V\big) \to C^*(\Gpd),
  $$
  from the Toeplitz algebra ${\cal T} \big(C(X),G,V\big)$ of our
interaction group \scite{\newsgrp}{5.1} to the groupoid C*-algebra
$C^*(\Gpd)$.  Since $C(X)\rtimes_V G$ is the quotient of the Toeplitz
algebra by the redundancy ideal \scite{\newsgrp}{6.2}, we need to show that 
$\widehat{\arep\times\lprep}$ vanishes on redundancies if we are to
reach our goal of 
obtaining  a homomorphism
  $$
  \arep\times\lprep : C(X)\rtimes_V G \to C^*(\Gpd).
  $$
  In the terminology of \scite{\newsgrp}{6.3} we must prove that
$(\arep,\lprep)$ is \stress{strongly covariant}.  The  following 
two Lemmas will be later used to further this goal.

But first let us introduce some notation.  By a \stress{word} in $G$
we shall mean any finite sequence $\a=(g_1,\ldots,g_k)$, where $g_i\in
G$.  If $\lprep$ is a partial representation of $G$ we shall let
  $$
  \prep\a = \prep{g_1}\ldots\prep{g_n}.
  $$

\state Lemma 
  \label CompleteToP
  Let $\lprep$ be a partial representation of $G$ in some
C*-algebra $B$ such that
$\prep n$ is an isometry for every $n\in P$.
  Then for every word $\a$ in $G$ there exists another word $\b$, such
that
  $$
  \prep\a\prep\b=\prep n,
  $$
  for some  $n\in P$.
  
  \proof Letting $\a=(g_1,\ldots,g_k)$ we shall prove the statement by
induction on $k$.
  If  $k=1$, in which case $\a=(g_1)$,  write $g_1=nm\inv$,
with $n,m\in P$, and let $\b=(m)$.  Then
  $$
  \prep\a\prep\b =
  \prep {g_1}\prep m =
  \prep {g_1}\prep m \prep m^*\prep m=
  \prep {{g_1}m} \prep m^*\prep m=
  \prep n,
  $$
  thus proving the result for $k=1$.

  If $k\geq 2$, let
$\a'=(g_2,\ldots,g_k)$ and use induction to get $\b'$ and $n'$ such
that $\prep {\a'}\prep {\b'}=\prep {n'}$.  Then
  $$
  \prep \a \prep {\b'} =
  \prep {g_1} \prep {\a'} \prep {\b'} =
  \prep {g_1} \prep {n'} = 
  \prep {g_1} \prep {n'}  \prep {n'}^* \prep {n'} =
  \prep {g_1n'} \prep {n'}^* \prep {n'} =
  \prep {g_1n'},
  $$
  and the conclusion follows from the case $k=1$.
  \proofend

The following idea has already been used in several occasions, e.g.~in 
\scite{\vershik}{8.6, 7.2}.

\state Lemma
  \label RightIdeal
  Given $n\in P$ there exists a finite set
$\{u_1,\ldots,u_n\}\subseteq C(X)$ such that
  $$
  \sum_{i=1}^n \arep(u_i) S_nS_n^*\arep(u_i) = 1.
  $$

\proof
  Let $\{V_i\}_{i=1}^n$ be a finite open cover of $X$ such that $\t_n$ is
injective on each $V_i$ and let $\{\phi_i\}_{i=1}^n$ be a partition of
unity subordinated to $\{V_i\}_{i=1}^n$.  Define $u_i(x) =
\phi_i(x)\cy(n,x)\mhalf$ and observe that for every $(x,g,y)\in\Gpd$
we have, using \lcite{\FormulaSSStar}, that
  $$
  \sum_{i=1}^n \arep(u_i) S_nS_n^*\arep(u_i)\calcat{(x,g,y)} =
  \sum_{i=1}^n u_i(x)u_i(y) S_nS_n^*\calcat{(x,g,y)}\$=
  \sum_{i=1}^n u_i(x)u_i(y) \cy(n,x) \half \cy(n,y)\half \bool{g=1}
\bool{\t_n(x)=\t_n(y)} \$=
  \sum_{i=1}^n u_i(x)^2 \cy(n,x) \bool{g=1} \bool{x=y} =
  \sum_{i=1}^n \phi_i(x)\bool{g=1} \bool{x=y} =
  \bool{g=1} \bool{x=y},
  $$
  proving the statement.
  \proofend

We now may prove the following crucial technical result:

  \def\M{{\cal M}}
  \def\rhoh{\widehat\rho}

\state Proposition
  \label TheCovarRepIsStronglyCovar
 The representation $(\arep,\lprep)$ of $\big(C(X),G,V\big)$ given by
\lcite{\TheCovarRep} is strongly covariant.

\proof Recall that to say that $(\arep,\lprep)$ is strongly covariant
is to say that  $\widehat{\arep\times\lprep}$ vanishes
on all redundancies \scite{\newsgrp}{6.3}.  In order to shorten our
notation we will write $\rhoh=\widehat{\arep\times\lprep}$.

Let $\a$ be a word in $G$ and let $k$ be an $\a$-redundancy.
Then $k\M_\a=0$, so
  $$
  0 = \rhoh(k\M_\a) =
  \rhoh(k) \rhoh(\M_\a) =
  \rhoh(k)\arep(C(X)) \prep \a \arep(C(X)),
  $$
  and in particular $\rhoh(k)\arep(C(X)) \prep \a=0$.
  Using Lemma \lcite{\CompleteToP} choose a word $\b$ such that
$\prep\a\prep\b=\prep n$, for some $n\in P$, so 
  $$
  0=
  \rhoh(k)\arep(C(X)) \prep \a\prep \b = 
  \rhoh(k)\arep(C(X)) \prep n.
  $$
  Recalling that $\prep n=S_n$ and using \lcite{\RightIdeal} we have
that
  $$
  \rhoh(k) =   \sum_{i=1}^n \rhoh(k)  \arep(u_i) S_nS_n^*\arep(u_i)  = 0.
  \proofend
  $$

Since $\widehat{\arep\times\lprep}$ vanishes on  redundancies it
factors through the quotient of ${\cal T}\big(C(X),G,V\big)$ by the
redundancy ideal and hence defines a *-homomorphism
  $$
  \arep\times\lprep: C(X)\rtimes_V G  \to C^*(\Gpd).
  $$

Our next major goal will be to prove that   $\arep\times\lprep$ is an
isomorphism.  The proof of injectivity  will be based on
\scite{\newsgrp}{10.6} so we are required to first  verify the
following:

  \state Proposition
  \label NonDeg
  The covariant representation $(\arep,\lprep)$ of $\big(C(X),G,V\big)$ given
by \lcite{\TheCovarRep} is non-degenerate.

  \proof
  Recall that to say that $(\arep,\lprep)$ is non-degenerate is to
say that the map
  $$
  a\in C(X) \mapsto \arep(a)\prep \a \in C^*(\Gpd)
  $$
  is injective for every word $\a$ in $G$.
  So suppose that $a\in C(X)$ is such that $\arep(a)\prep \a=0$.
Then, picking $\b$ and $n$ as in \lcite{\CompleteToP} we have
  $$
  0 = \arep(a) \prep \a \prep \b =
  \arep(a) \prep n =
  \arep(a) \S_n.
  $$
  Given $x\in X$ consider the element $(x,n,\t_n(x))\in\Gpd$.
Thinking of   $\arep(a) \S_n$ as a compactly supported function on
$\Gpd$ we compute
  $$
  0 = \arep(a)\prep n\calcat{(x,n,\t_n(x))} =
  a(x) \cy(n,x)\half,
  $$
  so it follows that $a=0$.
  \proofend

Notice that we have used that $\cy$ never vanishes in the proof above.
The following is the main result of this section:

  \state Theorem
  \label GroupoidModel
  In addition to  the hypotheses of \lcite{\StandingHyp} suppose that
$G$ is amenable and let
$\big(C(X),G,V\big)$ be the interaction group provided by
\lcite{\EndosGiveInteraction}, namely, if $g=n\inv m$, with $n,m\in
P$,
  $$
  V_g(f)\calcat y = \sum_{\theta_n(x)=y} \cy(n,x) f(\t_m(x))
  \for  f\in C(X)
  \for y\in X.
  $$
  Moreover consider the strongly covariant representation
$(\arep,\lprep)$ of $\big(C(X),G,V\big)$ in the groupoid C*-algebra of
$\Gpd$ given by \lcite{\TheCovarRep} and
\lcite{\TheCovarRepIsStronglyCovar}.  Then
  $$
  \arep\times\lprep : C(X)\rtimes_V G  \to C^*(\Gpd)
  $$
  is an isomorphism.

  \proof
  By \lcite{\NonDeg} and \scite{\newsgrp}{10.6} we have that
$\arep\times\lprep$ is injective on each $C_g$, where $\{C_g\}_{g\in
G}$ is the grading of $C(X)\rtimes_V G $ given by
\scite{\newsgrp}{7.2}.

Observe that $C^*(\Gpd)$ also  admits a grading $\{D_g\}_{g\in G}$
such that $\arep(f)\in D_1$, for all $f\in C(X)$, and $S_n\in D_n$, for
every $n\in P$.  If follows that $\arep\times\lprep$ is a graded
homomorphism in the sense that $(\arep\times\lprep)(C_g)\subseteq D_g$,
for every $g\in G$.  Let 
  $$
  \matrix{E_1:& C(X)\rtimes_V G & \to & C_1 \cr\cr
  E_2:& C^*(\Gpd) & \to & D_1}
  $$
  be the associated conditional expectations, so we have that
  $$
  E_2((\arep\times\lprep)(a)) = 
  (\arep\times\lprep)(E_1(a))
  \for a\in C(X)\rtimes_V G.
  $$
  If $a$ is such that $(\arep\times\lprep)(a)=0$, then 
  $$
  (\arep\times\lprep)(E_1(a^*a)) =   E_2((\arep\times\lprep)(a^*a)) =0.
  $$
  Since $E_1(a^*a)\in C_1$, and $\arep\times\lprep$ is injective on
$C_1$, as observed above, we have that $E_1(a^*a)=0$.  Given that $G$
is amenable, the Fell bundle $\{C_g\}_{g\in G}$ must also be amenable
\scite{\amena}{4.7}  and hence $E_1$ is faithful by
\scite{\amena}{2.12}.  So $a=0$.

In order to prove that $\arep\times\lprep$ is surjective observe that
the range of $\arep\times\lprep$ contains $\arep(f)$ and $S_n$, for
every $f\in C(X)$ and $n\in P$, so surjectivity  follows from
\lcite{\SpanGeneratesGpd}.
  \proofend

\section{Lattice-ordered semigroups}
  The reader might have a few examples in mind of semigroups acting on
compact spaces.  But before we can apply
\lcite{\EndosGiveInteraction}, or any of our  results based on it, 
we need to provide a normalized coherent cocycle, a task
which might not be entirely trivial.

In what follows we plan to show that coherent cocycles are often
present in a number of situations.
  With this goal in mind we will now study some elementary properties
of semigroups which will later play an important role in providing
applications of our results.

We will suppose throughout that $G$ is a group
and that $P$ is a subsemigroup of $G$ such that $P\cap P\inv=\{1\}$.
In most of our examples $G$ will be commutative but this does not seem
to be too significant for the general theory.  We shall therefore not
suppose that $G$ is commutative here.

  One may  define a left-invariant order on $G$ by saying that
  $$
  x\leq y \ \Longleftrightarrow \ x\inv y\in P.
  $$

\definition
  We shall say that the pair $(G,P)$ is a \stress{lattice-ordered
group} if, for every $x$ and $y$ in $G$, the set $\{x,y\}$ admits a
least upper bound $x\vee y$ and a greatest lower bound $x\wedge y$.

From now on we will fix a lattice-ordered group $(G,P)$.

\definition A \stress{mini-square} is by definition a quadruple of
elements $(s,t,u,v)\in P^4$ such that
  \izitem
  \zitem $su=tv$,
  \zitem $s \inf  t = 1$,
  \zitem $u\inv \sup v\inv = 1$.

  Observe that the last condition is equivalent to saying that 
  $u
  \mathop {\buildrel {\scriptscriptstyle r} \over \wedge}
  v = 1,
  $
  where 
  ``$\mathop {\buildrel {\scriptscriptstyle r} \over \wedge}$''
  denotes the greatest lower bound relative to the
\stress{right-invariant order} induced by $P$.
  Since it is a bit awkward to deal with two distinct order relations
at the same time we will make no further references to the
right-invariant order relation.

It is interesting to represent mini-squares by a diagram such as
  \medskip
  $$
  \matrix{ ^s\!\!\swarrow & \kern-6pt \searrow ^t\cr
           _u\!\!\searrow & \kern-6pt \swarrow _v} 
  $$
  \bigskip

\state Proposition 
  \label MNSTUV
  For every $m,n\in\G$ let
  $$s = (m\inf n)\inv m,\quad t = (m\inf n)\inv n,$$
  $$u = m\inv(m\sup n),   \quad v = n\inv(m\sup n).$$
  Then 
$(s,t,u,v)$  is a mini-square.

  \proof
  Since $(m\inf n)\leq m,n\leq (m\sup n)$ it is obvious that $s,t,u,v\in\P$.
  Obviously 
  $$ 
  su =  (m\inf n)\inv (m\sup n) = tv.
  $$
  In order to show that $s\inf t=1$, let $\g\in\G$ be such that 
  $ 
  \g\leq s,t.
  $
  Therefore 
  $$
  (m\inf n)\g \leq   (m\inf n)s = m
  \and
  (m\inf n)\g \leq   (m\inf n)t = n.
  $$
  It follows that
  $$
  (m\inf n)\g \leq     (m\inf n)
  $$
  and hence that $\g\leq e$. This shows that $s\inf t=1$.  Next we show
that $u\inv \sup v\inv = 1$.  For this suppose that
  $
  \g\geq u\inv,v\inv.
  $
  Then
  $$
  (m\sup n)\g\geq (m\sup n)u\inv =m
  \and
  (m\sup n)\g\geq (m\sup n)v\inv =n.
  $$
  Therefore
  $$
  (m\sup n)\g\geq   (m\sup n),
  $$
  whence $\g\geq e$.
  \proofend 

  \state Proposition 
  \label MiniSup
  Let $(s,t,u,v)$
  be a mini-square.  Then 
  $s\sup t = su = tv$.
  
  \proof
  By definition $su = tv$.    It is also clear that 
  $$
  su \geq s
  \and
  su = tv \geq t.
  $$
  Now suppose that $g\geq s,t$.  Then
  $$
  (su)\inv g \geq   (su)\inv s = u\inv
  \and
  (su)\inv g \geq   (su)\inv t = (tv)\inv t = v\inv.
  $$
  Therefore $$(su)\inv g \geq u\inv\sup v\inv = 1,$$ which implies that
$g\geq su$, as desired.
  \proofend

  \state Corollary
  \label CompletingMiniSquare
  Given $s,t\in\P$ such that $s\inf t=1$, there exists a unique pair
$(u,v)\in\P\times\P$ such that $(s,t,u,v)$ is a mini-square.

  \proof
  For existence apply \lcite{\MNSTUV} to the pair $(m,n):=(s,t)$.  For
uniqueness observe that by \lcite{\MiniSup} we have that
  $u=s\inv(s\sup t)$ and
  $v=t\inv(s\sup t)$.
  \proofend

Notice that 
if $G$ is commutative and $s,t\in G$ are such that $s\inf t=1$,
then
  $$
  \matrix{ ^s\!\!\swarrow & \kern-6pt \searrow ^t\cr
           _t\!\!\searrow & \kern-6pt \swarrow _s} 
  $$
  is clearly a mini-square.  By \lcite{\CompletingMiniSquare} these
are the only possible mini-squares.

\section{{\Admissible} semigroup actions and cocycles}
  The main goal of this section is to present sensible conditions on
semigroup actions and cocycles from which one may deduce coherence,
hence providing examples of interaction groups by
\lcite{\EndosGiveInteraction}.

We begin by fixing a lattice-ordered group
$(G,P)$.
  Observe that  \lcite{\EndosGiveInteraction.ii} is
automatically satisfied:

\state Proposition 
  \label PPinvInLattice
  If $(G,P)$ is a lattice-ordered group then $G=P\inv P=PP\inv.$

\proof
  Given $x\in G$, let $y=x\wedge e$.  Then $y\leq x$ and hence
$n:=y\inv x\in P$.  Observing that $y\inv\in P$,
because $y\leq e$, we have that $x=yn\in P\inv P$.  This shows that
$G=P\inv P$.  

Next let $n=x\vee e$.  Since $x\leq n$ we have that $m:=x\inv n\in P$.
So $x=n m\inv \in PP\inv$.
  \proofend

\definition
  \label AxiomX
  We shall say that a right action  $\t:P\to \End(X)$ is
\stress{\admissible}
  if, given any mini-square
  $$
  \matrix{ 
     ^s\!\!\swarrow & \kern-6pt \searrow ^t\cr
     _u\!\!\searrow & \kern-6pt \swarrow _v 
  } 
  $$
  and  $x,y\in X$ such that $\t_u(x)=\t_v(y)$, there exists a unique
$z\in X$ such that $\t_s(z)=x$, and $\t_t(z)=y$.

From now on we shall fix {\anadmissible} action $\t$ of $P$ on $X$.  Our
first task will be to describe sets of the form 
$\t_m\inv(p)\cap\t_n\inv(q).$  

\state Lemma
  \label InverseImagesIntersecting
  Given  $m,n\in P$ and $p,q\in X$, let  $(s,t,u,v)$ be the
mini-square given by \lcite{\MNSTUV} in terms of $m$ and $n$.
  \izitem
  \zitem If $\t_u(p)\neq\t_v(q)$, then $$\t_m\inv(p)\cap\t_n\inv(q)=\emptyset.$$
  \zitem If $\t_u(p)=\t_v(q)$, and $w$ is the unique element in $X$
such that $\t_s(w) =p$, and $\t_t(w)=q$, then 
  $$\t_m\inv(p)\cap\t_n\inv(q)=\t_{m\inf n}\inv(w).$$

  \proof
  (i) By contradiction let $x\in \t_m\inv(p)\cap\t_n\inv(q)$.  Then 
  $$
  \t_u(p) = 
  \t_u(\t_m(x)) =
  \t_{mu}(x) =
  \t_{nv}(x) =
  \t_v(\t_n(x)) =
  \t_v(q),
  $$
  contradicting the hypothesis.

  \medskip\noindent (ii)
  If $\t_{m\inf n}(x)=w$, then
  $$
  \t_m(x) =
  \t_{(m\inf n)s}(x) = 
  \t_s(\t_{m\inf n}(x)) = 
  \t_s(w) = p,
  $$
  and similarly
  $
  \t_n(x) =
  q.
  $
  Conversely, if 
  $\t_m(x) = p,$
  and 
  $\t_n(x) = q,$ set $w'=\t_{m\inf n}(x)$.  Then
  $$
  \t_s(w') = \t_s(\t_{m\inf n}(x)) = \t_{(m\inf n)s}(x) = \t_m(x) =p,
  $$
  and similarly
  $\t_t(w') = q.$
  By uniqueness we have that $w'=w$, so $\t_{m\inf n}(x)=w$, and hence
$x\in \t_{m\inf n}\inv(w).$
  \proofend

  \state Proposition
  \label ProductOfClasses
  Let  $\bar z\in X$, let $(s,t,u,v)$ 
  be a mini-square, and put $\bar x=\t_s(\bar z)$ and $\bar
y=\t_t(\bar z)$. Then the map
  $$
  \phi: \cl {s\sup t}{\bar z} \to  \cl u{\bar x}\times \cl v{\bar y},
  $$
  given by 
  $$
  \phi(z) = (\t_s(z),\t_t(z))
  $$
  is a bijection.

  \proof
  Given $z\in \cl {s\sup t}{\bar z}$, let $x=\t_s(z)$ and
$y=\t_t(z)$.  Then, recalling from \lcite{\MiniSup} 
that $s\sup t = su = tv$, notice that
  $$
  \t_u(x) = 
  \t_u(\t_s(z))=
\t_{su}(z) = \t_{s\sup t}(z) = \t_{s\sup t}(\bar z) =
\t_{su}(\bar z) =\t_u(\t_s(\bar z)) =\t_u(\bar x),
  $$
  proving that  $x\in \cl u{\bar x}$.  Similarly one shows that 
$y\in \cl v{\bar y}$, so $\phi(z)=(x,y)$ indeed lies in
$\cl u{\bar x}\times \cl v{\bar y}$.

  Let $(x,y)\in \cl u{\bar x}\times \cl v{\bar y}$.  Then 
  $$
  \t_u(x) = 
  \t_u(\bar x) = 
  \t_u(\t_s(\bar z)) = 
  \t_{su}(\bar z) = 
  \t_{tv}(\bar z) = 
  \t_v(\t_t(\bar z)) = 
  \t_v(\bar y) =
  \t_v(y).
  $$  
  By hypothesis there exists a unique $z\in X$ such that
$\phi(z)=(x,y)$.  We claim that $z$ must in fact lie in $\cl {s\sup
t}{\bar z}$.  To see this notice that
  $$
  \t_{s\sup t}(z) = 
  \t_{su}(z) =
  \t_u(\t_s(z)) = 
  \t_u(x) =
  \t_{su}(\bar z) =
  \t_{s\sup t}(\bar z).
  $$
  This shows that $\phi$ is bijective.
  \proofend

Let us now suppose we are given a normalized cocycle $\cy$. Given $x\in
X$ and $n \in P$ consider the restriction of  $\cy(n,\cdot)$ to  $\cl nx$.
By \lcite{\SumOne} we may view this as a probability distribution on
$\cl nx$.

Under the hypotheses of 
\lcite{\ProductOfClasses} observe that we then have two probability
distributions on $\cl {s\sup t}{\bar z}$, namely $\cy(s\sup t,\cdot)$ on
the one hand, and the product distribution of $\cy(u,\cdot)$ and
$\cy(v,\cdot)$ on the other.

Our next result touches upon the question as to whether these
probability distributions coincide.

\state Proposition
  \label SeveralAxioms
  Given a cocycle $\cy$, consider the following statements:
  \izitem
  \zitem for every $z\in X$ and every mini-square  $(s,t,u,v)$ one has that
  $\cy(s\sup t,z) =\hfill\break \cy(u,\t_s(z))\ \cy(v,\t_t(z)).$
  \zitem for every $z\in X$ and every mini-square  $(s,t,u,v)$ one has
that $\cy(t,z)=\cy(u,\t_s(z))$,
  \zitem for every $s,t\in\P$ with $s\inf t=1$, and for every
$z\in X$ one has that
  $\cy(s\sup t,z) = \cy(s,z)\ \cy(t,z)$.
  \medskip\noindent
  Then (i) $\Leftarrow$ (ii) $\Rightarrow$ (iii).  If $\cy$ never
vanishes then also  (i) $\Rightarrow$ (ii) $\Leftarrow$ (iii). 

  \proof
  \medskip\noindent (ii) $\Rightarrow$ (i).  We have
  $$
  \cy(s\sup t,z) = 
  \cy(tv,z) =
  \cy(t,z)\ \cy(v,\t_t(z)) =
  \cy(u,\t_s(z))\ \cy(v,\t_t(z)).
  $$

  \medskip\noindent (ii) $\Rightarrow$ (iii).  Given $(s,t)$ as in
(iii) pick $u$ and $v$ such that $(s,t,u,v)$ is a mini-square by \lcite{\CompletingMiniSquare}.
  We then have
  $$
  \cy(s\sup t,z) =
  \cy(su,z) =
  \cy(s,z)\ \cy(u,\t_s(z)) =
  \cy(s,z)\ \cy(t,z).
  $$

  \medskip\noindent  (i) $\Rightarrow$ (ii).  One has
  $$
  \cy(u,\t_s(z))\ \cy(v,\t_t(z)) =
  \cy(s\sup t,z) =
  \cy(tv,z) =
  \cy(t,z)\ \cy(v,\t_t(z)).
  $$   
  Since $\cy(v,\t_t(z))\neq 0$, we have that
  $\cy(u,\t_s(z)) = \cy(t,z).$

  \medskip\noindent (iii) $\Rightarrow$ (ii).  We have
  $$
  \cy(s,z)\ \cy(t,z) =
  \cy(s\sup t,z) = 
  \cy(su,z) = 
  \cy(s,z)\ \cy(u,\t_s(z)).
  $$
  Since   $\cy(s,z)\neq 0$ we conclude that 
  $\cy(t,z) = \cy(u,\t_s(z))$.
  \proofend

  We therefore make the following:

  \definition
  A cocycle $\cy$ is said to be \stress{\admissible}
  if \lcite{\SeveralAxioms.ii} and hence also
\lcite{\SeveralAxioms.i \& iii} holds.

The reason why we are interested in  {\admissible} cocycles is given below.

\state Proposition
  \label CoherentFromAdmissible
  If $\t$ and $\cy$ are {\admissible} then $\cy$ is coherent.

\proof  We first claim that if  $\t_{m\sup n}(z) \neq\t_{m\sup n}(x)$,
then both sides of \lcite{\DefineCoherent} vanish.  By symmetry it is
enough to show that this is so for the right-hand side of 
\lcite{\DefineCoherent}.

  Let $(s,t,u,v)$ be the mini-square obtained from \lcite{\MNSTUV}
from $m$ and $n$, and set $p=\t_m(z)$, and $q=\t_n(x)$.
Then 
  $$
  \t_u(p)=
  \t_u(\t_m(z)) = 
  \t_{mu}(z) =
  \t_{m\sup n}(z) \neq
  \t_{m\sup n}(x) =
  \t_{nv}(x) =
  \t_v(\t_n(x)) =
  \t_v(q).
  $$  
  That is, $\t_u(p) \neq \t_v(q)$, and hence by
\lcite{\InverseImagesIntersecting} we have that
  $$
  \cl nx\cap \cl mz = \t_n\inv(q) \cap \t_m\inv(p) = \emptyset.
  $$
  The right-hand side of  \lcite{\DefineCoherent} then vanishes because
$\Cl nmxz=\emptyset$. 

  Suppose now that $\t_{m\sup n}(z) = \t_{m\sup n}(x)$.  Letting
$(s,t,u,v)$, $p$ and $q$ be as above we than conclude similarly that
$\t_u(p)=\t_v(q)$.  Again  by
\lcite{\InverseImagesIntersecting} we have that
  $$ 
  \cl nx\cap \cl mz = \t_n\inv(q) \cap \t_m\inv(p) =
  \t_{m\inf n}\inv(w),
  $$
  where $w$ is the unique element in $X$ such that $\t_s(w) =p$, and
$\t_t(w)=q$.  Therefore the right-hand side of \lcite{\DefineCoherent}
satisfies
  $$
  \sum_{y\in \cl mz\cap \cl nx} \cy(m,y)\, \cy(n,x) = 
  \sum_{y\in \t_{m\inf n}\inv(w)} \cy((m\inf n)s,y)\, \cy(n,x) \$=   
  \sum_{y\in \t_{m\inf n}\inv(w)} \cy(m\inf n,y) \,\cy(s,\t_{m\inf n}(y))
\, \cy(n,x) \$=
  \left(\sum_{y\in \t_{m\inf n}\inv(w)} \cy(m\inf n,y)\right) \cy(s,w)
\, \cy(n,x) =
  \cy(s,w) \, \cy(n,x) = 
  (\star).
  $$
  Observe that by symmetry   $(t,s,v,u)$ is a mini-square as well.
Thus, applying \lcite{\SeveralAxioms.ii}, we have that
  $$
  \cy(s,w) \={(\SeveralAxioms.ii)}
  \cy(v,\t_t(w)) =
  \cy(v,\t_n(x)).
  $$
  Therefore 
  $$
  (\star) = \cy(v,\t_n(x))\, \cy(n,x) = \cy(nv,x) = \cy(m\sup n,x).
  $$
  Having arrived at an expression which is symmetric in the variables
$m$ and $n$, the proof is complete.
  \proofend

  We thus arrive at the following  important result:
  
\state Corollary
  \label MainCorollary
  Let $(G,P)$ be a lattice-ordered group, let $\t$ be {\anadmissible}
right action of $P$ on a compact space $X$, and let $\cy$ be
{\anadmissible} normalized cocycle.  Then there exists a unique
interaction group $V=\{V_g\}_{g\in G}$ on $C(X)$, such that, if
$g=n\inv m$, with $n,m\in P$, then
  $$
  V_g(f)\calcat y = \sum_{\theta_n(x)=y} \cy(n,x) f(\t_m(x))
  \for  f\in C(X)
  \for y\in X.
  $$

  \proof
  Follows immediately from the above result plus
\lcite{\PPinvInLattice} and \lcite{\EndosGiveInteraction}.
  \proofend

\section{Example: single endomorphism}
  \label SingleEndoSection
  Beginning with this section we shall give several examples of
{\admissible}, therefore coherent, cocycles for which one may apply
\lcite{\GroupoidModel}.

  Let $\G=\Z$ and let $\P=\N$.  It is obvious that $(\Z,\N)$ is a
lattice-ordered group.
  Given any endomorphism $T$ of $X$ define an action $\t$ of $\N$ on $X$ by
  $$
  \t:n\in\N\mapsto T^n\in\End(X).
  $$

The only mini-squares in $\N$ are of the form 
  $$
  \matrix{ ^n\!\!\swarrow & \kern-6pt \searrow ^0\cr
           _0\!\!\searrow & \kern-6pt \swarrow _n} 
  $$
  where $n$ $\in\N$,
  and its  reflections across the vertical axis,
  and hence it is obvious that $\t$ is {\admissible}.

  For each $x\in X$, let 
  $$
  \cy(0,x)=1
  \and
  \cy(1,x) = {1\over |\cl 1x|}.
  $$
  For all $n\geq 1$, define 
  $\cy(n+1,x)$ recursively by 
  $$
  \cy(n+1,x)=\cy(n,x)\ \cy(1,\t_n(x)).
  $$
  It is then easy to prove by induction that $\cy$ is a normalized
{\admissible} cocycle which does not vanish anywhere.

We should observe that, in the case of the present example, Theorem
\lcite{\GroupoidModel} was already proved in \scite{\vershik}{9.1}.

\section{Example: star-commuting endomorphisms} 
  \label StarComutSect
  This example is inspired by \cite{\ArzRena}.
\definition
  \label DefStarCommuted
  Let $S$ and $T$ be commuting maps on a set $X$.  We shall say that
the pair $(S,T)$ \stress{star-commutes} \cite{\ArzRena} if for every
$x,y\in X$ such that $T(x)=S(y)$, there exists a unique $z\in X$ such
that $S(z)=x$, and $T(z)=y$.

  Consider $\G=\Z\times\Z$ and let $\P=\N\times\N$.  It is
obvious that $(\Z\times\Z,\N\times\N)$ is a lattice-ordered group with
  $$
  (n_1,m_1)\sup(n_2,m_2) = (\max\{n_1,n_2\},\max\{m_1,m_2\})
  $$
  and
  $$
  (n_1,m_1)\inf(n_2,m_2) = (\min\{n_1,n_2\},\min\{m_1,m_2\}).
  $$
  Observe that all mini-squares are either given by 
  $$
  \matrix{ ^{(n,m)}\!\!\swarrow & \kern-6pt \searrow ^{(0,0)}\hfill\cr
           \hfill_{(0,0)}\!\!\searrow & \kern-6pt \swarrow _{(n,m)}} 
  \hbox {\quad or \quad } 
  \matrix{ ^{(n,0)}\!\!\swarrow & \kern-6pt \searrow ^{(0,m)}\cr
           _{(0,m)}\!\!\searrow & \kern-6pt \swarrow _{(n,0)}} 
  $$
  where $n,m\in\N$, or the reflections of the above across a vertical axis.

Let $S$ and $T$ be commuting endomorphisms of a compact Hausdorff
space $X$ and define
  $$
  \t: (n,m)\in\N\times\N \mapsto S^{n}T^{m}\in\End(X).
  $$
  It is clear that $\t$ is an action of $\N\times\N$ on $X$.
  Observe that \lcite{\DefStarCommuted} is then equivalent to saying
that the condition of \lcite{\AxiomX} holds for the special
mini-square
  $$
  \matrix{ ^{(1,0)}\!\!\swarrow & \kern-6pt \searrow ^{(0,1)}\cr
           _{(0,1)}\!\!\searrow & \kern-6pt \swarrow _{(1,0)}} 
  $$

\state Proposition
  \label StarImplyAdmiss
  If $(S,T)$ star-commutes then $\t$ is {\admissible}.

  \proof
  Leaving aside the trivial mini-squares (i.e., those involving the
trivial group element) and taking advantage of vertical symmetry, we
consider only mini-squares of the form
  $$
  \matrix{ ^{(n,0)}\!\!\swarrow & \kern-6pt \searrow ^{(0,m)}\cr
           _{(0,m)}\!\!\searrow & \kern-6pt \swarrow _{(n,0)}} 
  $$
  where $n,m\geq1$.
  Our task therefore consists in showing the following:

  \medskip
  ``If \ $T^m(x)=S^n(y)$,
then there exists a unique $z\in X$ such that $S^n(z)=x$, and $T^m(z)=y$.''

  \medskip
We shall prove this by induction on $n+m$.  Given that $n,m\geq1$, we
see that the lowest possible value for  $n+m$ is 2, in which case the
claim follows directly from the hypothesis.

  Assuming $n+m\geq3$, suppose without lack of generality that $n\geq2$.
  Then, letting $y'=S(y)$, notice that 
  $T^m(x)=S^{n-1}(y')$,
  so we have by the induction hypothesis that there exists a unique $w\in X$ such that
  $$
  S^{n-1}(w)=x\and T^m(w)=y'.
  \eqno{(\dagger)}
  $$

  Next, observing that $T^m(w)=y'=S(y)$, and that $1+m<n+m$, we see, again by induction,
that there exists $z\in X$ such that 
  $$
  S(z)=w\and T^m(z)=y.
  \eqno{(\ddagger)}
  $$
  Clearly
  $$
  S^n(z)=
  S^{n-1}(S(z)) =
  S^{n-1}(w) = x.
  $$
  Thus $z$ solves our existence question.  As for uniqueness assume
that $z'$ is such that $S^n(z')=x$, and $T^m(z')=y$.  Setting
$w'=S(z')$, we have
  $$
  S^{n-1}(w') =   S^{n-1}(S(z')) = S^n(z') = x
  $$ 
  and
  $$
  T^m(w')=   T^m(S(z'))= S(T^m(z'))= S(y)=y'.
  $$
  Comparing the above with \lcite{$\dagger$}, and using the uniqueness
part of the induction hypothesis, we see that necessarily $w'=w$.  Observing that
  $$
  S(z')=w
  \and
  T^m(z')=y,
  $$
  we then conclude, in view of \lcite{$\ddagger$} and the induction
hypothesis, that $z'=z$.
  \proofend

  We next construct an {\admissible} hence coherent cocycle.  Let
  $$
  \cy_S,\cy_T:\N\times X \to \R_+,
  $$
  be each given as in section \lcite{\SingleEndoSection} relatively to 
$S$ and $T$, respectively.  

  \state Lemma
  \label InvariantCocycle
  If $(S,T)$ star-commutes then
  for every $n\in\N$ and every $x\in X$ one has
  $$
  \cy_S(n,x) =   \cy_S(n,T(x)) 
  \and
  \cy_T(n,x) =   \cy_T(n,S(x)).
  $$
  
  \proof
  We prove the first assertion only, doing  so by induction on
$n$.  Obviously it holds for $n=0$.  Speaking of the case $n=1$, let
$x\in X$.   Adopting the notation 
  $$
  \cl Sx = \{x'\in X:S(x')=S(x)\},
  $$
  consider the map
  $$
  \phi: x'\in \cl Sx \mapsto T(x') \in \cl S{T(x)}.
  $$
  That $T(x')$ in fact belongs to $\cl S{T(x)}$ follows from
  $$
  S(T(x')) = 
  T(S(x')) = 
  T(S(x)) = 
  S(T(x)).
  $$
  We claim that this map is bijective.  In order to show that it is
one-to-one
let $x',x''\in \cl Sx$ be
such that $T(x')=T(x'')$.
  Consider the diagram
  $$
  \qquad\qquad
  \matrix{
  & \square \cr
  \hfill^S \swarrow && \searrow^T\hfill \cr
  \vrule height16pt depth 10pt width0pt
  S(x)\quad && T(x')=T(x'')\cr
  \hfill_T \searrow && \swarrow_S\hfill \cr
  &  \kern-8pt S(T(x)) \kern-8pt 
  }
  $$
  Since both $x'$ and $x''$ fit in the box we conclude that $x'=x''$,
by the fact that $(S,T)$ star-commutes.

  In order to show surjectivity let $y\in \cl S{T(x)}$ and consider the diagram
  $$
  \matrix{
  & \square \cr
  \hfill^S \swarrow && \searrow^T\hfill \cr
  \vrule height16pt depth 10pt width0pt
  S(x)\quad && \kern18pt y\cr
  \hfill_T \searrow && \swarrow_S\hfill \cr
  &  \kern-8pt S(T(x)) \kern-8pt 
  }
  $$
  By hypothesis there exists a $x'$ which fits in
the box.  Obviously such a $x'$ lies in $\cl Sx$ and $T(x')=y$.  This
proves that our map is in fact bijective.

  As a consequence we have that
  $|\cl Sx| = |\cl S{T(x)}|$, so that
  $$
  \cy_S(1,x) =
  {1\over |\cl Sx|} =
  {1\over |\cl S{T(x)}|} =
  \cy_S(1,T(x)).
  $$
  Given $n\geq1$ we have by induction that
  $$
  \cy_S(n+1,T(x)) =
  \cy_S(n,T(x))\ \cy_S(1,S^n(T(x))) =
  \cy_S(n,x)\ \cy_S(1,T(S^n(x))) \$=
  \cy_S(n,x)\ \cy_S(1,S^n(x)) =
  \cy_S(n+1,x).
  \proofend
  $$  

We should observe that the cocycles $\cy_S$ and $\cy_T$ given by
section \lcite{\SingleEndoSection} are certainly not unique.
Nevertheless, if a different choice of $\cy_S$ and $\cy_T$ was made in
the above proof it is not clear that we could carry it on.  It is
therefore interesting to pinpoint exactly to what extent are those
cocycles special.

  \state Theorem 
  \label StarCommutExistCocycle
  If $(S,T)$ star-commutes then the map
  $
  \cy_T:(\N\times\N)\times X \to \R_+,
  $
  defined by 
  $$
  \cy((n,m),x)=\cy_S(n,x)\ \cy_T(m,x)
  $$
  is a normalized {\admissible} cocycle which does not vanish anywhere.
  
  \proof
  It is obvious that  $\cy((0,0),x)=1$, for all $x$, so let us 
check the cocycle identity \lcite{\TransferLaw}.  Given
$n=(n_1,n_2)$ and $m=(m_1,m_2)$ in $\N\times\N$ we have
  $$
  \cy(n+m,x)= 
  \cy((n_1+m_1,n_2+m_2),x)=
  \cy_S(n_1+m_1,x)\ \cy_T(n_2+m_2,x) \$=
  \cy_S(n_1,x)\ \cy_S(m_1,S^{n_1}(x))\
    \cy_T(n_2,x)\ \cy_T(m_2,T^{n_2}(x)) =$$$$ \={(\InvariantCocycle)}
  \cy_S(n_1,x)\ \cy_S(m_1,T^{n_2}S^{n_1}(x))\
    \cy_T(n_2,x)\ \cy_T(m_2,S^{n_1}T^{n_2}(x)) \$=
  \cy((n_1,n_2),x)\ \cy((m_1,m_2),S^{n_1}T^{n_2}(x))=
  \cy(n,x)\ \cy(m,\t_n(x)).
  $$

We next prove that $\cy$ is {\admissible}.  With respect to mini-squares
of the form 
  $$
  \matrix{ ^{(n,0)}\!\!\swarrow & \kern-6pt \searrow ^{(0,m)}\cr
           _{(0,m)}\!\!\searrow & \kern-6pt \swarrow _{(n,0)}} 
  $$
  condition \lcite{\SeveralAxioms.ii} takes  the form
  $\cy_T(m,z)=\cy_T(m,S^n(z))$, while for mini-squares of the form 
  $$
  \matrix{ ^{(0,m)}\!\!\swarrow & \kern-6pt \searrow ^{(n,0)}\cr
           _{(n,0)}\!\!\searrow & \kern-6pt \swarrow _{(0,m)}} 
  $$
  it becomes
  $\cy_S(n,z)=\cy_S(n,T^m(z))$.  That $\cy$ is {\admissible} then follows
immediately from Lemma \lcite{\InvariantCocycle}.

  It remains to prove \lcite{\SumOne}.  So let $(n,m)\in \N\times\N$
and  $y\in  X$.  Then
  $$
  \sum_{\t_{(n,m)}(x)=y}\cy((n,m),x) =
  \sum_{S^n(T^m(x))=y}\cy_S(n,x)\,\cy_T(m,x) \$=
  \sum_{S^n(w)=y}\ \sum_{T^m(x)=w}\cy_S(n,T^m(x))\,\cy_T(m,x) =
  \sum_{S^n(w)=y}\cy_S(n,w) \sum_{T^m(x)=w}\cy_T(m,x) = 1.
  \proofend
  $$

  \section{Ledrappier's shift}
  We now wish to give an interesting concrete example of a pair of
star-commuting endomorphisms of Bernoulli's space.
  Let us first recall the construction of Ledrappier's dynamical
system \cite{\Ledrappier}, which is also discussed in \cite{\Deaconu}.
  Let
  $$
  K= \{0,1\}^{\N^2}
  $$ 
  have the product topology 
  and let $X$ be the subset of $K$ formed by the elements $x\in K$
such that 
  $$
  x_{n,m} +   x_{n+1,m} +   x_{n,m+1}  = 0 
  \for (n,m)\in\N^2,
  \eqno{(\dagger)}
  $$
  where addition is performed modulo 2.
  Clearly $X$ is a compact subspace of $K$ which  is invariant under both the
horizontal and vertical shifts, namely the transformations
  $$
  H,V : K \to K
  $$
  given by 
  $$
  H(x)_{n,m} = x_{n+1,m}
  \quad \hbox{and} \quad
  V(x)_{n,m} = x_{n,m+1}.
  $$

  Our short term goal is to show that the restrictions of $H$ and $V$
to $X$ give a pair of star-commuting endomorphisms, for which we may
then apply the results of section \lcite{\StarComutSect}.

  Let us view the elements of $K$ as possible ways of arranging zeros
and ones along the vertices of the lattice $\N^2$.  It is easy to see
that $X$ consists precisely of such arrangements which are entirely
made of the following four patterns
  \bigskip
  $$
  \matrix{
    \cdot & \cdot & \cdot & \cdot \cr
    \cdot & 0 & \cdot & \cdot \cr
    \cdot & 0 & 0 & \cdot \cr
    \cdot & \cdot & \cdot & \cdot 
}
\qquad \qquad\qquad \qquad
  \matrix{
    \cdot & \cdot & \cdot & \cdot \cr
    \cdot & 1 & \cdot & \cdot \cr
    \cdot & 0 & 1 & \cdot \cr
    \cdot & \cdot & \cdot & \cdot 
}
  $$
\bigskip
  $$
  \matrix{
    \cdot & \cdot & \cdot & \cdot \cr
    \cdot & 1 & \cdot & \cdot \cr
    \cdot & 1 & 0 & \cdot \cr
    \cdot & \cdot & \cdot & \cdot 
}
\qquad \qquad\qquad \qquad
  \matrix{
    \cdot & \cdot & \cdot & \cdot \cr
    \cdot & 0 & \cdot & \cdot \cr
    \cdot & 1 & 1 & \cdot \cr
    \cdot & \cdot & \cdot & \cdot 
}
  $$
  \bigskip
  If on $(\dagger)$ we perform the change of variables
$(p,q)=(n,m+1)$, we may rewrite it as
  $$
  x_{p,q} =  -x_{p,q-1} -x_{p+1,q-1}
  \for p\in\N \for q\geq1,
  $$
  which says that if $x\in X$, then each  entry of $x$ is determined
by the entry immediately below it and the one to the right of that.
In particular each   row of $x$  is determined by the row below it and
hence the first row determines everything.
  In other words the map
  $$
  \{x_{n,m}\}_{n,m\in\N} \in X \longmapsto 
  \{x_{n,0}\}_{n\in\N} \in \Omega
  $$
  is a homeomorphism from $X$ onto Bernoulli's space $\Omega=\{0,1\}^\N$.  It
is easy to see that, under this homeomorphism, the horizontal shift
$H$ identifies  with the usual Bernoulli shift $S$, given by
  $$
  S : \{\w_n\}_{n\in\N} \in \Omega \longmapsto
\{\w_{n+1}\}_{n\in\N}\in\Omega.  \qquad
  $$
  Almost as easily one checks that the vertical shift $V$ becomes the
map
  $$ 
  T : \{\w_n\}_{n\in\N} \in \Omega \longmapsto
\{\w_n+\w_{n+1}\}_{n\in\N}\in\Omega.
  $$
  Observe that $\Omega$ is an abelian group under coordinatewise
addition modulo 2.  With this group structure one has that $S$ and $T$
are group homomorphism and  $T=id+S.$  
  
  \def\Ker{{\rm Ker}}%
  \state Lemma
  Let $G$ be an abelian group, let $\phi:G\to G$ be a surjective group
homomorphism, and let $\psi=\phi+id$.  Then $(\phi,\psi)$ is
star-commuting.
  
  \proof
  It is obvious that $\phi$ and $\psi$ commute.
  Next  suppose that $\phi(x) = \psi(y)$.  Choose $z_1\in G$ such that
$\phi(z_1)=y$, and notice that
  $
  k:= x-y-z_1 \in \Ker(\phi).
  $
  Indeed
  $$
  \phi(k) =
  \phi(x)-\phi(y)-\phi(z_1) =
  \phi(x)-\phi(y)-y =
  \phi(x) - \psi(y) =
  0.
  $$
  Setting $z=z_1+k$, one obviously has that $\phi(z) = \phi(z_1+k) =
\phi(z_1) = y$, and moreover
  $$
  \psi(z) =
  \phi(z) + z =
  y + z_1 + x-y-z_1 =
  x.
  $$
  In order to show uniqueness, suppose that $z'\in G$ is such that
$\phi(z')=y$ and $\psi(z')=x$.  Then
  $$
  z' =
  \psi(z') - \phi(z') =
  x-y=
  \psi(z) - \phi(z) =
  z.
  \proofend
  $$

As an immediate corollary we conclude that $(S,T)$ is star-commuting.
It is well known that $S$ is an endomorphism of $\Omega$.   It can
be shown that $T$ is conjugated to $S$ and hence $T$ is also an
endomorphism.  We therefore  find ourselves precisely in the situation
of section \lcite{\StarComutSect}, as desired.

  \section{Example: an action of the multiplicative integers} 
  Consider the multiplicative group of strictly positive rational numbers
$\G=\Q_+^\times$, and the sub-semigroup $\P=\N\setminus\{0\}$.
  It is
clear that $(\G,\P)$ is a lattice-ordered group with
  $$
  n\sup m = \gcd(n,m)
  $$
  (greatest common divisor)
  and
  $$
  n\inf m = \lcm(n,m)
  $$
  (least common multiple).

  Observe that, in the present context, the mini-squares  are given by 
  $
  \matrix{ ^n\!\!\swarrow & \kern-6pt \searrow ^m\cr
           _m\!\!\searrow & \kern-6pt \swarrow _n} 
  $
  where $n$ and $m$ are relatively prime.
  Consider the action $\t$ of\/ $\P$ on the unit circle $S^1$ given by
  $$
  \t_n(x) = x^n
  \for x\in S^1.
  $$

  \state Proposition
  The action $\t$ defined above  is {\admissible}.

  \proof Let $(n,m,m,n)$ be a mini-square such as the one above and
suppose that $x,y\in S^1$ are such that $x^m=y^n$.  We need to find
$z\in S^1$ such that
  $$
  z^n=x \and z^m=y.
  $$
  Let $w$ be any $n^{th}$ root of $x$, so that $w^n=x$, and notice that
  $$
  (w^{-m} y)^n = 
  (w^n)^{-m} y^n =
  x^{-m} y^n = 1.
  $$
  Thus $\lambda:= w^{-m} y$ is an $n^{th}$ root of 1.   
  Since $\gcd(n,m)=1$, the mapping $\rho\mapsto \rho^m$ is a bijection
on the set of all $n^{th}$ roots of 1.  Thus we may choose an $n^{th}$
root $\rho$ of 1 such that $\rho^m=\lambda$.  Setting  $z = \rho w$ we have
  $$  
  z^n = \rho^n w^n = x 
  \and
  z^m =
  \rho^mw^m =
  \rho^m \lambda\inv y =
  y,
  $$
  thus taking care of existence.  If $z_1$ and $z_2$ are elements  of $S^1$
such that
  $$
  z_1^n=x=z_2^n \and z_1^m=y=z_2^m,
  $$
  write $z_j=e^{2\pi\theta_j}$, for $j=1,2$, and observe that
  $$
  e^{2\pi n\theta_1}= e^{2\pi n\theta_2} \implies p:=
n(\theta_2-\theta_1)\in\Z,
  $$
  and similarly $q:=m(\theta_2-\theta_1)\in\Z$.  Thus
  $$
  {p\over n}=
  (\theta_2-\theta_1) = 
  {q\over m},
  $$
  so that 
  $nq=mp$.  Since $\gcd(n,m)=1$, we have that $n$ divides
$p$, so that 
  $
  \theta_2-\theta_1\in\Z
  $
  and hence $z_1=z_2$.
  \proofend

  There is a very elementary cocycle we can define in the present
context:

  \state Proposition Setting
  $\cy(n,x) = 1/n$, we have that 
  $\cy$ is a normalized {\admissible} cocycle which does not vanish
anywhere.

  \proof
  Left to the reader.
  \proofend

\section{Example: polymorphisms} 
  Let $X$ be a compact space. 
  According to \cite{\ArzRena} a polymorphism 
  $$
  X\buildrel T \over \leftarrow\Sigma\buildrel S \over \to X
  $$
consists of a pair of surjective local homeomorphisms $S,T:\Sigma\to
X$, where $\Sigma$ is a compact space.  The idea is to think of
``$ST\inv$" as a multivalued map from $X$ to $X$.

One special case of interest \cite{\ArzRena} is when $\Sigma=X$ and
$ST=TS$.  In this case we may generalize the usual notion of
transformation groupoid by introducing the following:
  \def\HGpd{{\cal H}}
  $$
  \HGpd =
  \big\{ (x,k,y)\in X\times \Z\times X: \exists\, n,m\in\N,\ k=n-m,\
S^nT^m(x) = S^mT^n(y)\big\}.
  $$

Note that the definition of the groupoid of a polymorphism given in
\cite{\ArzRena} is not correct and should be replaced by this one.

  \state Proposition
  $\HGpd$ is a groupoid under the operations
  $$
  (x,k,y) (y,\ell,z) = (x,k+\ell, z)
  \and 
  (x,k,y)\inv = (y,-k,x).
  $$

  \proof Given $(x,k,y)$ and $(y,\ell,z)$ in $\HGpd$, let $k=n_1-m_1$,
and $\ell = n_2-m_2$, such that 
  $$
  S^{n_1}T^{m_1}(x) = S^{m_1}T^{n_1}(y)
  \and
  S^{n_2}T^{m_2}(y) = S^{m_2}T^{n_2}(z).
  $$
  Then
  $$
  S^{n_1+n_2}T^{m_1+m_2}(x) = 
  S^{n_2}T^{m_2} S^{n_1}T^{m_1}(x) =
  S^{n_2}T^{m_2} S^{m_1}T^{n_1}(y) \$=
  S^{m_1}T^{n_1} S^{n_2}T^{m_2}(y) =
  S^{m_1}T^{n_1} S^{m_2}T^{n_2}(z) = 
  S^{m_1+m_2}T^{n_1+n_2}(z),
  $$
  proving that 
  $$
  \HGpd \ni (x,(n_1+n_2)-(m_1+m_2),z) = (x,k+\ell,z).
  $$
  That $(y,-k,x)\in\HGpd$ is obvious.
  We leave it for the reader to verify the remaining points.
  \proofend

Alternatively consider the action
$\t$ of $\N\times\N$ on $X$ defined near the beginning of section
\lcite{\StarComutSect} by
  $$
  \t_{(n,m)} = S^nT^m
  \for (n,m)\in\N\times\N,
  $$
  and let $\Gpd$ be the groupoid defined in terms
of $\t$ as in \lcite{\IntroduceGroupoid}.

One should not expect $\HGpd$ and $\Gpd$ to coincide since the former
should be thought of as the transformation groupoid associated to the
dynamical system generated by the ``multivalued" map $ST\inv$, while
the latter is generated by $S$ and $T$.  But $\HGpd$ may be viewed as
a subgroupoid of $\Gpd$ in the following sense:

  \def\d{d}
  \state Proposition 
  Let $\d:\Gpd\to\Z$ be the groupoid homomorphism given by
  $$
  d(x,(n,m),y) = n+m.
  $$
  Then the map
  $$
  \phi:(x,k,y)\in\HGpd \ \longmapsto\ (x,(k,-k),y)\in\Gpd
  $$
  defines an isomorphism from $\HGpd$ to the kernel of $d$.  
  
  \proof
  Let us begin by verifying that $\phi(x,k,y)$ does indeed belong to
$\Gpd$,  for every $(x,k,y)\in\HGpd$.  
Write $k=n-m$, in such a way that $S^nT^m(x) = S^mT^n(y)$.  Since this
is equivalent to
  $\t_{(n,m)}(x) =   \t_{(m,n)}(y)$,
  we have that
  $$
  \Gpd\ni
  (x,(n,m)-(m,n),y)=
  (x,(k,-k),y).
  $$

  The only non obvious remaining point is perhaps to prove that the
kernel of $\d$ is contained in the range of $\phi$.  In order to prove
this let $(x,(n,m),y)\in{\rm Ker}(d)$, so that $m=-n$.  By definition
there exists $(p,q),(r,s)\in\N\times\N$ such that
  $$
  (n,m)=(n,-n) =
(p,q)-(r,s)
  \and 
  \t_{(p,q)}(x) = \t_{(r,s)}(y).
  $$
  In other words $n=p-r=s-q$, and 
 $S^pT^q(x)=S^rT^s(y)$.  For any $(k,l)\in\N\times\N$ notice that
  $$
  S^{k+p}T^{l+q}(x)=
  S^kT^lS^pT^q(x) =
  S^kT^lS^rT^s(y) =
  S^{k+r}T^{l+s}(y).
  $$
  Choose  $l\in\N$ big enough for $l+s-p\geq0$, and set $k=l+s-p$.  Then
obviously $k+p=l+s$, while
  $$
  l+q =
  l+s-p+r =
  k+r,
  $$
  thus proving that 
  $$
  \HGpd \ni
  (x,(k+p)-(l+q),y) = (x,n,y).
  $$
  Moreover
  $$
  \phi(x,n,y) = (x,(n,-n),y) = (x,(n,m),y).
  \proofend
  $$

Suppose from now on that $(S,T)$ star-commutes.  Then $\t$ is an
admissible action by \lcite{\StarImplyAdmiss}, and hence $\Gpd$ is a
locally compact \'etale groupoid by \lcite{\Topology}.  It is clear
that $\phi(\HGpd)$ is an open subgroupoid of $\Gpd$ and hence $\HGpd$
may be turned into a locally compact \'etale groupoid with the induced
topology.  Moreover let $\cy$ be a never vanishing normalized coherent
cocycle for $\t$, such as that given by
\lcite{\StarCommutExistCocycle}.  One may then consider the
interaction group $(C(X),\Z\times\Z,V)$ given by
\lcite{\EndosGiveInteraction}.  Define
  $$
  \mu: n\in\Z\ \longrightarrow\ (n,-n)\in\Z\times\Z,
  $$
  and let
  $$
  W_n = V_{\mu(n)}
  \for n\in\Z.
  $$
  An explicit formula for $W_n$ is easy to write down:
  $$
  W_n(f)\calcat y = \left\{\matrix{
    \ds\sum_{T^n(x)=y} \hfill \cy(n,x) f(S^n(x)) \hfill, & \hbox{ if }
n\geq0, \cr\cr\cr
    \ds\sum_{S^{-n}(x)=y} \cy({-n},x) f(T^{-n}(x)) \hfill, & \hbox{ if
} n<0,} \right.
  $$
  for every $f\in C(X)$, and $y\in X$.  
  For example, $W_1(f)\calcat y$ is the weighted average of $f(z)$ as
$z$ runs in the finite set $S(T\inv(\{y\}))$.

  $W_1(f)$ is  therefore  the single valued
function that best mimics the composition of $f$ with the multivalued
function $ST\inv$.  The latter being precisely what one wants
to capture by considering the polymorphism
  $X\buildrel T \over \leftarrow\Sigma\buildrel S \over \to X$.

  It is then obvious that $(C(X),\Z,W)$ is an interaction group.  It
is interesting to notice that, contrary to what happens to $V$, and
unless either $T$ or $S$ are invertible maps, there
is no $n\in\Z$ for which $W_n$ is an endomorphism of $C(X)$.  

\state Theorem
  There is a natural isomorphism
  $$
  \phi: C(X)\rtimes_W \Z \to C^*(\HGpd).
  $$

  \proof
  The proof follows essentially the same method used to prove
\lcite{\GroupoidModel}, so we restrict ourselves here to a brief
outline.

  First observe that $\HGpd$ is open in $\Gpd$, so $C^*(\HGpd)$ is
naturally a subalgebra of $C^*(\Gpd)$.  Considering the covariant
representation $(\arep,\lprep)$ given by \lcite{\TheCovarRep} it is
immediate to verify that $(\arep,\lprep\circ\mu)$ is a covariant
representation of $(C(X),\Z,W)$, which actually takes values within
$C^*(\HGpd)$.  That it is strongly covariant is an immediate
consequence of $(\arep,\lprep)$ possessing this property.  We
therefore obtain a *-homomorphism
  $$
  \arep\times(\lprep\circ\mu) : C(X)\rtimes_W \Z \to C^*(\HGpd),
  $$
  which is injective by the  same reasons used in
\lcite{\TheCovarRep}.  Surjectivity also follows as before.
  \proofend

\section{A counter-example}
  In this section we want to show an example of a semigroup action
which does not admit a never vanishing coherent cocycle.

Let us begin by introducing some notation: given a compact space $X$
and a function $S:X\to X$ define an equivalence relation $R_S$ on $X$
by
  $$
  (x,y)\in R_S  \iff S(x)=S(y).
  $$
  Recall that if $R_1$ and $R_2$ are relations on $X$ one defines
the \stress{composition} of $R_1$ and $R_2$ to be the relation
  $$
  R_1\circ R_2 = 
  \{(x,z)\in X\times X: \exists\, y\in X,\ (x,y)\in R_1,\ (y,z)\in R_2\}.
  $$
  We will  say that $R_1$ and $R_2$ commute if
  $R_1\circ R_2 = R_2\circ R_1$.

  \state Proposition
  Let $G$ be a group, $P$ be a subsemigroup of $G$, and
$\t:P\to\End(X)$ be a semigroup action.  If  $\cy$ is  a never
vanishing coherent cocycle for  $\t$, then $R_{\t_n}$ and $R_{\t_m}$
commute for every $n,m\in P$.

  \proof
  Observe that $(x,z)\in R_{\t_n}\circ R_{\t_m}$ if and only if $\Cl
{n}{m} xz\neq\emptyset$ (see \lcite{\DefineBigCl} for a definition
of $\Cl {n}{m} xz$).
  Given $(x,z)\in R_{\t_n}\circ R_{\t_m}$ we then have that
  $\W{m}{\Cl {n}{m} xz}>0$, because $\cy$ is never zero.
Since $\cy$ is coherent we deduce from \lcite{\DefineCoherent} that
  $\W{n}{\Cl {m}{n} xz}>0$, and hence that $\Cl {m}{n}
xz\neq\emptyset$, from were we have that $(x,z)\in R_{\t_m}\circ
R_{\t_n}$.  This shows that
  $$
  R_{\t_n}\circ R_{\t_m} \subseteq R_{\t_m}\circ R_{\t_n}.
  $$
  The converse inclusion follows similarly.
    \proofend

In order to exhibit a semigroup action which does not admit a never
vanishing coherent cocycle it is therefore enough to provide one for
which the conclusion of the above result fails.

Given commuting endomorphisms $S,T\in\End(X)$, define a semigroup
action of $\N\times\N$ on $X$ by
  $$
  \theta: (n,m)\in\N\times\N \mapsto S^{n}T^{m}\in\End(X),
  $$
  as in section \lcite{\StarComutSect}.
  Our plan is thus to provide endomorphisms $S$ and $T$ for which $R_S$
and $R_T$ do not commute.

  For every  $p\in\N$, let $\O_p:= \{0,1\}^p$ and let $\O=\O_\infty
=\{0,1\}^\N$.  With the product topology $\O$ is a compact space also
known as Bernoulli's space.
  Let $S$ be the shift on $\O$, namely the transformation
  $$
  S:(x_0,x_1,x_2,\ldots)\in\O\ \longmapsto\ (x_1,x_2,x_3,\ldots) \in\O.
  $$
  It is well known that $S$ is an endomorphisms of $\O$.

  We next wish
to describe a class of transformations $T:\O\to\O$ which commute with
$S$.  For this let $p\in\N$ and choose any subset $D\subseteq \O_p$,
henceforth referred to as the \stress{dictionary}.
  For $x=(x_0,x_1,x_2,\ldots)\in\O$ define 
  $$
  T(x)_k = \left\{\matrix{
    1, & \hbox{ if } (x_{k},x_{k+1},\ldots,x_{k+p-1}) \in D, \cr \cr
    0, &  \hbox{ otherwise.}\hfill
    }\right.
  $$

In other words, to compute $T(x)$ one slides a \stress{window} of
width $p$ along $x$ and checks whether or not the word seen through
the window belongs to the dictionary.  Recording the answers as
a sequence of ones (when the word belongs to $D$) and zeros
(otherwise), we obtain $T(x)$.

\bigskip 
\begingroup 
\noindent \hfill\beginpicture
\setcoordinatesystem units <0.0015truecm, 0.0015truecm> point at 0 0

\put{$(\ x_0,x_1,x_2,\ldots,x_{p-1},x_p,
  \ldots, x_k,x_{k+1},\ldots,x_{k+p-1},x_{k+p}, \ldots\ )$} at 5000 400
\setdashpattern <2pt,1.5pt> 
\def\square#1#2#3#4{\plot #1 #3 #2 #3 #2 #4 #1 #4 #1 #3 /}
\square{1300}{3680}{190}{600}
\put{{\fiverm initial window position}}  at 2500 700 
\square{4720}{7280}{190}{600}
\setsolid
\arrow <0.15cm> [0.35,0.75] from 5600 700 to 6400 700

\endpicture \hfill\null 
\endgroup

\bigskip
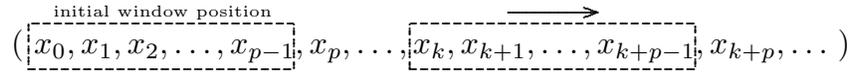
\centerline{\eightrm  The sliding window method.} 
\bigskip\medskip

It is easy to see that $T$ is a continuous mapping which commutes with
$S$.  Such transformations are sometimes called \stress{cellular
automata}.  By a well known result due to Hedlund \cite{\Hedlund}, any
continuous map $T:\O\to\O$ which commutes with $S$ has the above form
for a suitable dictionary $D$.

We now wish to describe a specific class of dictionaries for which the
associated cellular automaton is a surjective local homeomorphism.

\definition
  We shall say that a given subset $D\subseteq\O_p$ is
\stress{progressive} if for every $\b\in\O_{p-1}$ there exists a
unique $\varepsilon\in\{0,1\}$ such that $\b\varepsilon\in D$
(concatenation). 

Notice that when $D$ is progressive and we are about to choose the
last component of a word which already has $p-1$ components, say
  $$
  \b = (x_0,\ldots,x_{p-2},\ ?\ ),
  $$
  we always  have the option of forming a word which belongs to $D$ or not. 

  \state Theorem
  \label Surjective
  Let $T$ be a cellular automaton associated to a progressive
dictionary $D$.  Then $T$ is a surjective local homeomorphism, that
is, $T\in\End(\O)$.

  \proof
  For the duration of this proof we shall find it useful to extend
$T$ to finite words, obtaining, for every $m\geq p$, a map
  $$
  T_m: \O_m \to \O_{m-p+1},
  $$
  defined by  the above sliding window method.

  Given $y\in\O$, let $\b=(x_0,\ldots,x_{p-2})$ be any word in
$\O_{p-1}$.  Since $D$ is progressive there exists a unique
$x_{p-1}\in \{0,1\}$ such that 
  $$
  \b_0:= (x_0,\ldots,x_{p-2},x_{p-1})\in D \iff y_0=1.
  $$
  One then has that 
  $
  T_{p}(\b_0)
  $
  is the single-letter word $(y_0)$.
  Again because $D$ is progressive there exists a unique
$x_p\in\{0,1\}$ such that
  $$
  (x_1,\ldots,x_{p-1},x_p) \in D \iff
  y_1=1.
  $$
  Thus
  $$
  \b_1:  = (x_0,x_1,\ldots,x_{p-1},x_p) \in\O_{p+1},
  $$
  and $T_{p+1}(\b_1) = (y_0,y_1)$.
  By continuing in this fashion one builds an infinite word
$x\in\Omega$ such that $T(x)=y$, hence proving that $T$ is surjective.

Observe that, after the initial  choice of $\b$, each $x_k$ (for
$k\geq p-1$) was uniquely determined.  Therefore the restriction of
$T$ to 
  $$
  V_\b = \{x\in\O : x_i=\b_i,\ \forall i=0,\ldots,p-2\},
  $$
  gives a bijection onto $\O$.
  Since $V_\b$ is compact we conclude that $T|_{V_\b}$ is a homeomorphism.
Moreover the collection $\{V_\b\}_{\b\in\O_{p-1}}$ is an open cover of
$\O$, and hence $T$ is a local homeomorphism.
  \proofend

\sysstate{Example}{\rm}{
  \label MyNicestExample
  Consider the progressive set 
  $$
  D=\{000,100,010,111\}\subseteq\O^3,
  $$
  and let $T$ be the cellular automaton associated to $D$.  Then
$(S,T)$ is a pair of endomorphisms of $\O$ such that $R_S$ and $R_{T}$
do not commute.
  }

  \proof
  Let $x,y,z\in \O$ be given by
  $$
  x = 0\ 1\ 1\ 1\ 1\ 1\ldots
  \ ,\quad  
  y = 1\ 1\ 1\ 1\ 1\ 1\ldots
  \ \and
  z = 0\ 0\ 0\ 0\ 0\ 0\ldots
  $$
  Then
  $$
  S(x) = S(y) = 
  T(y) = T(z) = 1\ 1\ 1\ 1\ 1\ldots
  $$
  so $(x,y)\in R_S$, and $(y,z)\in R_{T}$, whence 
  $(x,z)\in R_S\circ R_{T}$.
  We claim that $(x,z)\notin R_{T}\circ R_S$.  By contradiction
suppose otherwise so there exists $y'$ such that
  $$
  T(x) = T(y') \and
  S(y') = S(z).
  $$
  The second equation allows for only two choices for $y'$, namely
  $$
  y'_1 = 0\ 0\ 0\ 0\ 0\ 0\ldots
  \and
  y'_2 = 1\ 0\ 0\ 0\ 0\ 0\ldots
  $$
  none of which satisfy $T(x) = T(y')$. \proofend

\state Corollary
  Let $S$ and $T$ be as in \lcite{\MyNicestExample} and define a
semigroup action of\/ $\N\times\N$ on $\O$ by
  $$
  \theta: (n,m)\in\N\times\N \mapsto S^{n}T^{m}\in\End(\O).
  $$
  Then there does not exist a never vanishing coherent cocycle for
$\t$.

  It seems that this example is as odds with the last Proposition in
\cite{\Deaconu}.

  \endgroup

  \begingroup
  \bigskip\bigskip \bigskip  \font \sc = cmcsc8 \sc
  \parskip = -1pt

  Departamento de Matem\'atica \hfill D\'epartment de Math\'ematiques

  Universidade Federal de Santa Catarina \hfill Universit\'e d'Orl\'eans

  88040-900 -- Florian\'opolis -- Brasil \hfill 45067 Orl\'eans, France

  \eightrm exel@\kern1pt mtm.ufsc.br \hfill Jean.Renault@univ-orleans.fr

  \endgroup
  \end